\documentclass[twocolumn]{autart}

\usepackage {amsmath,times}
\usepackage {amssymb}
\usepackage {graphicx}
\usepackage{amsrefs}
\usepackage{color}
\usepackage{tikz}
\usetikzlibrary{arrows}
\usepackage{subcaption}
\usepackage{hyperref}
\usepackage{wrapfig}
\usetikzlibrary{positioning}

\newtheorem{lemma}[thm]{Lemma}
\newtheorem{proposition}[thm]{Proposition}
\newtheorem{corollary}[thm]{Corollary}
\newtheorem{definition}[thm]{Definition}
\newtheorem{remark}[thm]{Remark}

\newcommand{\SX}{\mathbb{X}}

\newcommand{\mj}[1]{{ \color{red} \sf $\clubsuit$ Mauricio J.: [#1]}}

\begin{document}

\begin{frontmatter}

\title{On reachability of Markov chains:  A long-run average approach}

\author[D. Avila]{Daniel \'Avila}	\ead{daniel.avila@uclouvain.be},
\author[M. Junca]{Mauricio Junca}\ead{mj.junca20@uniandes.edu.co}

\address[D. Avila]{Center for Operations Research and Econometrics, Universit\'e Catholique de Louvain, Louvain-la-Neuve, Belgium.}
\address[M. Junca]{Department of Mathematics, Universidad de los Andes, Bogot\'a, Colombia.}


\begin{abstract}
We consider a Markov control model in discrete time with countable both state space and action space. Using the value function of a suitable long-run average reward problem, we study various reachability/controllability problems. First, we characterize the domain of attraction and escape set of the system, and a generalization called $p$-domain of attraction, using the aforementioned value function. Next, we solve the problem of  maximizing the probability of reaching a set $A$ while avoiding a set $B$. Finally, we consider a constrained version of the previous problem where we ask for the probability of reaching the set $B$ to be bounded. In the finite case, we use linear programming formulations to solve these problems. Finally, we apply our results to a example of an object that navigates under stochastic influence.
\end{abstract}

\end{frontmatter}



\section{Introduction}

Markov Decision Processes (MDP) provide a mathematical framework for modeling decision making problems where outcomes are uncertain. Formally, these are discrete time Markov control stochastic processes. The most common and studied setting are stochastic systems over a discrete state space with discrete action space, see for example ~\cites{puterman2014markov, bertsekas1995dynamic}, but general spaces are also part of the literature, see ~\cites{bertsekas1978stochastic,hernandez1996discrete, hernandez2012further}. We will focus on the former case. Applications of MDPs range from inventory control and investment planning to economics and behavioral ecology.

Based on the problem there are different planning horizons involved. In finite horizon problems one is interested in the evolution of the process upon a time $T$, while in infinite horizon problems the interest is in the long time behavior of the system. In every case, the objective is to maximize the expected cumulative reward obtained over time, where such reward depends on the state of the system and the action taken at each time (sometimes it also includes the state at the next period).

In this work we will be interested in using the MDP framework to solve some problems concerning the controllability/reachability of a controlled Markov chain. As we will see, such problems do not only evaluate the state of the system at each time, but on the whole evolution of the process. The first problem we aim to solve is to characterize the domain of attraction and escape set of a closed set $A$. The first one refers to the initial states for which there exist a control that takes the system to $A$, while the escape set are the initial states such that no control can take the system to $A$. A somehow related problem was studied in \cite{arevalo,arvelo2883control}, where the idea was to use entropy methods to maximize the number of recurrent states. Let $\SX$ be the state space of the Markov model, so we define these sets as follows.
\begin{definition}\label{DomAtt&EscSet}
Given a set $A\subset\SX$, let
\begin{align*}
\Lambda_A=\left\lbrace x\in \mathbb{X}\Big|  \liminf_{t\rightarrow\infty } P_x^{\pi}(X_t\in A )>0 \text{ for some policy }\pi \right\rbrace,\\
\Gamma_A =\left\lbrace x\in \mathbb{X}\Big| \liminf_{t\rightarrow\infty } P_x^{\pi}(X_t\in A)=0 \text{ for all policies }\pi \right\rbrace.
\end{align*}
\end{definition}
The domain of attraction  appears in the context of deterministic differential equations when describing the initial states under which the system will approach to a stable point. Such technique is commonly referred in the literature as Zubov's method, see ~\cites{zubov1964methods,hahn1967stability}. It allows to characterize the domain of attraction  and the escape set in terms of an appropriate value function, which is the solution of a differential equation. In ~\cite{camilli2008control}, Zubov's method is generalized for deterministic controlled systems, and in~\cite{camilli2006zubov} it is further generalized for stochastic differential equations. For the present work we took as guide the constructions made in this last work. Inspired by the literature of stochastic target problem, see \cites{Bouchard,Soner}, we study the $p$-domain of attraction for any $p\in[0,1]$, defined as follows.
\begin{definition}\label{pDomAtt}
Given a set $A\subset\SX$ and $p\in(0,1]$, let
\begin{align*}
\Lambda_{A,p}=\left\lbrace x\in \mathbb{X}\mbox{ } \Big|  \liminf_{t\rightarrow\infty } P_x^{\pi}(X_t\in A )\geq p \mbox{ for some policy }\pi \right\rbrace.
\end{align*}
\end{definition}
A similar problem in the context of stochastic hybrid systems and finite horizon is studied in \cite{abate}. 

The next problem consists on finding a control policy that maximize the probability of reaching some set $A$ while avoiding some set $B$. Namely, let $\tau_A$ and $\tau_{B}$ be the hitting times of $A\subset\SX$ and $B\subset\SX$, respectively. Consider an initial distribution $\nu$ over the state space. Our main objective will be to find a control policy $\pi$ that solves the problem
\begin{equation}\tag{P1}\label{P1}
  \max\limits_{\pi}P^{\pi}_{\nu}(\tau_A<\tau_B, \tau_A<\infty).
\end{equation}
Note that the set $B$ acts as a cemetery set since the evolution of the controlled Markov process is meaningless after this set is reached. This problem is studied in \cites{chatterjee2011maximizing}, where the authors consider general state and actions spaces but assume that $\tau_A\wedge\tau_B<\infty$ a.s for every policy $\pi$, hence they use a total reward approach. Finite horizon related problems for hybrid systems and a continuous time version for controlled diffusions can be found in \cite{summers2010verification} and \cite{Esfahani}, respectively.

Finally, we consider a constrained version of the previous problem as follows:
\begin{align}\label{P2}\tag{P2}
  &\max\limits_{\pi}  P^{\pi}_{\nu}(\tau_A<\infty) \\\notag
  &\mbox{s.t. }  P^{\pi}_{\nu}(\tau_B<\infty)\leq \epsilon.
\end{align}
In this case the set $B$ is no longer a cemetery set, making the evolution of the controlled process different, depending on whether the set $B$ has been reached or not. This fact suggests that this cannot be a Markovian problem. To the best of our knowledge, such problem, or any similar, has not been studied in the literature.

The contributions of this work are the following:
\begin{enumerate}
\item[(i)] We consider reachability problems in infinite horizon and relate them with long-run average reward problems in the context of MDPs. The main result in this direction is Theorem \ref{eqAv-Ht} that calculates the probability of reaching a closed set in finite time in terms of such reward problems. Then, we use this result to characterize $p$-domains of attraction (Corollary \ref{DoAvf}) and solutions of reach-avoid problems (Theorem \ref{ThMarkov}). An important feature of our results is that it includes the case multi-chain models. To the best of our knowledge, this is the first time that long-run average reward problems are used to characterize reachability problems. 
\item[(ii)] We consider a reachability problem with a constraint in the probability of hitting a given set of states. In general, stochastic control problems with probabilistic constrains are hard to solve since in this case Dynamic Programming Principle is unclear. Using a state space augmentation technique, we are able to formulate the problem in terms of long-run average reward problems (Corollary \ref{corP2}). In the case of finite state and action spaces we use linear programming duality to solve the problem (Theorem \ref{finite_case}).
\end{enumerate}

The paper is organized as follows: Section \ref{prelim} consists of two parts. In subsection \ref{MDP} we state the framework and known results for Markov Decision Processes, in particular linear programming formulations for long-run average reward problems. Subsection \ref{closed} presents the result that relates the probability of reaching closed set in finite time with the MDP reward problem. Sections \ref{domain}, \ref{max} and \ref{maxcon} presents the results for $p$-domains of attraction, Problem \eqref{P1} and Problem \eqref{P2}, respectively. In Section \ref{lpFinite}, we consider finite state and action spaces and formulate Problem \eqref{P2} as a linear program. Finally, in Sections \ref{example} and \ref{future} we present a numerical example to illustrate our findings and some future research directions.

\section{Markov control model and closed sets}\label{prelim}

In this section we establish the framework of Markov Decision Processes (MDP) and some results about closed sets that allow to formulate properly the problems described above.

\subsection{Markov Decision Processes}\label{MDP}

We will describe now the decision model, see \cite{hernandez1996discrete,puterman2014markov} for details. A Markov control model is a tuple $(\mathbb{X},\mathbb{U}, \{\mathbb{U}(x) | x\in\mathbb{X}\}, Q , r )$, where:
\begin{enumerate}
\item $\mathbb{X},\mathbb{U}$ are sets corresponding to the state space and actions respectively; in our work $\mathbb{X}$ and $\mathbb{U}$ are countable. 
\item For each $x\in\mathbb{X}$ there is a set $\mathbb{U}(x)\subset\mathbb{U}$ corresponding to the feasible actions to state $x$. The set of feasible states and actions is denoted as $\mathbb{K}$, that is,
\[ \mathbb{K}=\{ (x,u) | x\in\mathbb{X}, u\in\mathbb{U}(x)\}.
\]
\item $Q$ is a stochastic kernel on $\mathbb{X}$ given $\mathbb{K}$, that is, for each $(x,u)\in\mathbb{K}$, $Q(\cdot|x,u)$ is a probability measure on $\mathbb{X}$ and for each $B\subset\SX$, the function $Q( B|\cdot)$ is measurable on $\mathbb{K}$. 
\item $r$ is a function $r:\mathbb{K}\rightarrow\mathbb{R}$, called the reward function, which we assume bounded.
\end{enumerate}

A control policy is a sequence $\pi=\{\mu_t\}_{t\geq0}$ of stochastic kernels on $\mathbb{U}$ given the set of admissible histories up to time $t$, $H_t:=\mathbb{K}^t\times\SX$, such that for all $t\geq0$ and $h_t\in H_t$
$$\mu_t\left(\mathbb{U}(x_t)|h_t=(x_0,u_0,\ldots,x_{t-1},u_{t-1},x_t)\right)=1.$$
We denote by $\Pi$ the set of control policies.  

Consider the measurable space $(\Omega,\mathcal{F})$, where $\Omega:= (\mathbb{X}\times\mathbb{U})^{\mathbb{N}}$ and $\mathcal{F}$ is the product sigma algebra. Given a policy $\pi\in\Pi$ and a distribution $\nu$ over $\mathbb{X}$ there exist a unique probability measure $P_{\nu}^{\pi}$ on $(\Omega,\mathcal{F})$ and a $\SX\times\mathbb{U}$-valued stochastic process $(\{(X_t,U_t)\}_{t\geq0},P_{\nu}^{\pi})$ such that for each $B\subset\mathbb{X}$ and $V\subset\mathbb{U}$, (see \cite{hernandez1996discrete})
\begin{itemize}
\item $P_{\nu}^{\pi}(\mathbb{K}^{\mathbb{N}})=1$
\item $P_{\nu}^{\pi}(X_0\in B)=\nu(B)$
\item $P_{\nu}^{\pi}(U_t\in V|h_t)=\mu_t(V|h_t)$
\item $P_{\nu}^{\pi}(X_{t+1}\in B|h_t,U_t=u)=Q(B|x_t,u).$
\end{itemize}
We denote by $E_{\nu}^{\pi}$ the expected value with respect to $P_{\nu}^{\pi}$. When $\nu=\delta_x$ for $x\in\SX$ we use the notation $P_x^{\pi}$ and $E_x^{\pi}$. In general the process $\{X_t\}_{t\geq0}$ is non-Markovian, so, in order to make the process Markovian we also consider Markovian control policies. In this case we have a sequence $\pi=\{\mu_t\}_{t\geq0}$, where $\mu_t$ is a stochastic kernel on $\mathbb{U}$ given only $\SX$. Denote by $\Pi_M\subset\Pi$ the set of Markovian control policies. In this case the process $\{X_t\}_{t\geq0}$ is Markovian and
$$P_{x,x'}^{\pi} :=P^{\pi}_{\nu}(X_{t+1}= x'| X_t=x) = \sum_{u\in\mathbb{U}(x)} \mu_t(u|x)Q(x'|x,u).$$
If $\mu_t$ is the same for all $t\geq0$ the Markovian policy is called stationary and in this case $\{X_t\}_{t\geq0}$ is a time-homogeneous Markov process. Finally, if $\pi$ is stationary and $\mu(u|x)\in\{0,1\}$ for all $u\in\mathbb{U}(x)$ and $x\in\mathbb{X}$, then we say that $\pi$ is a stationary deterministic policy.

We will denote by $\overline{x}^t$ the history up to time $t$ of the process $\{X_t\}$, and similarly for the process $\{U_t\}$. Hence, admissible histories can be written as $h_t=(\overline{x}^{t-1},\overline{u}^{t-1},x_t)$. We also say that $\overline{x}^t\notin B$ if none of its components belong to $B$.

\subsubsection{Long-run average reward problems}\label{rewprob}

Given a Markov control model there are three infinite horizon control problems studied in the literature: Discounted, total and long-run average reward problems. We will focus on the latter case. Given an initial state $x\in\SX$, we want to find a policy $\pi\in\Pi$ that maximizes the reward function
\[
  v^{\pi}(x) :=\limsup_{N\rightarrow\infty}\frac{1}{N} \mathbb{E}_x^{\pi}\Big[ \sum_{t=0}^{N-1}r(X_t,U_{t}) \Big].
\] 
The first important result about these problems is that for any $\pi\in\Pi$ there exists $\pi'\in\Pi_M$ such that 
\begin{equation}\label{eqmarginals}
P_{\nu}^{\pi}(X_t=x',U_t=u)=P_{\nu}^{\pi'}(X_t=x',U_t=u)
\end{equation}
for any $t\geq0$, $x'\in\mathbb{X}$ and $u\in\mathbb{U}(x')$, which implies that the reward functions are the same, see \cite{puterman2014markov}. Therefore, optimal policies can always be found among Markovian policies.

Given $x\in\SX$, let $V(x)=\sup\limits_{\pi\in\Pi}v^{\pi}(x)$. The existence of an optimal policy relies on the so-called optimality equations for the multi-chain model: For all $x\in\mathbb{X}$
\begin{align}\label{MC1}\tag{MC1}
 & \max_{u\in\mathbb{U}(x)} \Big\{ \sum_{x'\in\mathbb{X}} P_{x,x'}^u v(x') - v(x)\Big\} =0,\\\label{MC2}\tag{MC2}
 &\max_{u\in \mathbb{U}(x)} \Big\{ r(x,u) - v(x) + \sum_{x'\in\mathbb{X}} P_{x,x'}^u h(x') - h(x) \Big\} = 0,
\end{align}
where $P_{x,x'}^u:=Q(x'|x,u)$.

\begin{remark}
Multi-chain models are those where there exists a stationary Markovian policy for which the induced Markov chain has at least two recurrent classes. Uni-chain models are easier since there is just one recurrent class and the function $V$ is constant, so there is no need to introduce the first equation. Here, we consider the multi-chain case since we will deal with closed sets.
\end{remark}
We have the following theorem that relates a solution of equations above with the average reward function, for a proof see Section 9.1 of \cite{puterman2014markov}.
\begin{thm}\label{AveragePrimal}
\begin{enumerate}
\item Suppose that the optimality equations have a solution $(v,h)$, then, $v=V$. 
\item Suppose that the state space $\mathbb{X}$ is finite and $\mathbb{U}(x)$ is finite for all $x\in\mathbb{X}$. Then, the optimality equations have a solution.
\item Suppose that the state space $\mathbb{X}$ is finite and $\mathbb{U}(x)$ is finite for all $x\in\mathbb{X}$. Then, there exist a stationary deterministic policy $\pi$ such that $v^{\pi}=V$. Moreover, such optimal policy can be defined as $\pi=\{\mu\}$ where
\[
     \mu(u|x)\in \arg\max_{u\in\mathbb{U}(x)}\{ r(x,u) + \sum_{x'\in\mathbb{X}}P_{x,x'}^uh(x') \}
\]
with the additional condition that $P^{\pi}v=v$, where $P^{\pi}$ is the transition matrix induced by $\pi$.
\end{enumerate}
\end{thm}

When the state space $\mathbb{X}$ is finite and $\mathbb{U}(x)$ is finite for all $x\in\mathbb{X}$, the optimality equations can also be solved via linear programming. Consider a vector $\nu\in \mathbb{R}_{\geq 0}^{|\mathbb{X}|}$. The linear program is as follows, see \cite{puterman2014markov} for further details:
\begin{align}\label{AvPrimal}\tag{AvP}
  \min & \sum_{x\in\mathbb{X}} \nu(x) v(x) \\\notag
  \mbox{s.t. }& \\\notag
 & v(x) \geq \sum_{x'\in\mathbb{X}} P_{x,x'}^u v(x'),\\\notag
 & v(x)  \geq  r(x,u)+ \sum_{x'\in\mathbb{X}} P_{x,x'}^u h(x') - h(x),\\
 & \forall x\in\mathbb{X},u\in\mathbb{U}(x).
\end{align}
Also, its dual program is given by
\begin{align}\label{AvDual}\tag{AvD}
  &\max \sum_{x\in\mathbb{X}} \sum_{u\in\mathbb{U}(x)}r(x,u)\alpha(x,u)  \\\notag
  &\mbox{s.t.} \\\notag
 & \sum_{u\in\mathbb{U}(x)} \alpha(x,u) - \sum_{x'\in\mathbb{X}\atop u\in\mathbb{U}(x')} P_{x',x}^u \alpha(x',u)  = 0,\\\notag
 & \sum_{u\in\mathbb{U}(x)} \alpha(x,u)  +   \beta(x,u) - \sum_{x'\in\mathbb{X}\atop u\in\mathbb{U}(x')} P_{x',x}^u \beta(x',u)  = \nu(x), \\\notag
 & \alpha(x,u)\geq 0 ,\beta(x,u)\geq 0,\quad\forall x\in\mathbb{X},u\in\mathbb{U}(x).
\end{align}
The next theorem relates the above linear programs with the multi-chain optimality equations (see a proof in Section 9.3 of \cite{puterman2014markov}).
\begin{thm}\label{AverageDual}
\begin{enumerate}
\item There exist optimal solutions $v^*$ and $(\alpha^*,\beta^*)$ of the primal and dual linear program, respectively.
\item Define the stationary policy $\pi=\{\mu\}$ as
\[
\mu(u^+|x)=
\begin{cases}
\dfrac{\alpha^* (x,u^+)}{\sum\limits_{u\in\mathbb{U}(x)}\alpha^*(x,u)},  \mbox{ if } \sum\limits_{u\in\mathbb{U}(x)}\alpha^*(x,u)>0 \\
\\
\dfrac{\beta^* (x,u^+)}{\sum\limits_{u\in\mathbb{U}(x)}\beta^*(x,u)}, \mbox{ otherwise.  }  \\
\end{cases}
\]
Then, $v^{\pi}=V=v^*$. 
\end{enumerate}
\end{thm}

\subsection{Closed sets}\label{closed}

Closed subsets of the state space are essential for the correct formulation of the problems described in the introduction. We start with the definition of a closed set under a policy $\pi$ which says that once the process hits the set it stays in the set afterwards.
\begin{definition}
Given a policy $\pi\in\Pi$ and a set $A\subset \mathbb{X}$, we say $A$ is closed under $\pi$ if given that $x_s\in A$ for some $s\leq t$, then
$$P_{\nu}^{\pi}(X_{t+1}\notin A|h_{t})= \sum_{u_{t}\in \mathbb{U}(x_{t})} Q(A^c|x_{t},u_{t})\mu_t(u_{t}|h_{t} )= 0.$$
\end{definition}
\begin{remark}\label{ClosedQ}
Note that if $Q(x'|x,u)=0$ for all $x\in A$, $x'\notin A$ and $u\in\mathbb{U}(x)$, then $A$ is a closed set under any policy $\pi\in\Pi$.
\end{remark}

Now, given a set $A\subset\mathbb{X}$, the random variable $\tau_A= \inf \{t\geq 0 \mbox{ }|\mbox{ } X_t\in A\}$ is called the hitting time of $A$. We have the following result about closed sets and its proof is included in Appendix \ref{appA}.

\begin{proposition}\label{HtA}
Assume $A$ is a closed set under $\pi\in\Pi$. Then, for any $x\in \mathbb{X}$
\[\lim_{t\rightarrow\infty } P^{\pi}_x(X_t\in A )= P^{\pi}_x(\tau_A<\infty ).
\]
\end{proposition}

The importance of the previous proposition is that it allows to express the probability of an event that depends on the joint distribution of the process, in terms of probabilities of events that only depend on the marginal distributions. Therefore, we have the following theorem, which is the main result of this section. Consider the Markov control model described in Subsection \ref{MDP} with reward function $r=\mathbf{1}_A$, the indicator function of set $A$. Hence, for $x\in\SX$ and $\pi\in\Pi$ the associated long-run average reward is given by
\begin{equation}\label{rewardA}
 v^{\pi}(x) :=\limsup_{N\rightarrow\infty}\frac{1}{N}\mathbb{E}_x^{\pi}\left[ \sum_{t=0}^{N-1}\mathbf{1}_{A}(X_t)\right].   
\end{equation}

\begin{thm}\label{eqAv-Ht}
Let $x\in\mathbb{X}$ and assume that $A$ is a closed set under $\pi\in\Pi$. Then, $v^{\pi}(x) =P_x^{\pi}(\tau_A<\infty)$.
\end{thm}
\begin{pf} For any $N\in\mathbb{N}$ we have that
\begin{align*}
      \frac{1}{N}\mathbb{E}_x^{\pi}\Big[ \sum_{t=0}^{N-1}\mathbf{1}_{A}(X_t)\Big]  &= \frac{1}{N}\sum_{t=0}^{N-1} P_x^{\pi}(X_t\in A).
\end{align*}   
Taking $\limsup$ on both sides we get
\begin{align*}
v^{\pi}(x) &= \limsup_{N\rightarrow \infty}\frac{1}{N}\sum_{t=0}^{N-1}P_x^{\pi}(X_t\in A).
\end{align*}
Now, recall that given a sequence $\{s_t\}_{t\geq0}$ such that the limit $L$ exists, then the Ces\`aro limit also exists and it is equal to $L$, that is, 
\[
  \lim_{N\rightarrow\infty}\frac{s_0+\ldots+s_{N-1}}{N} =L.
\] 
Since $A$ is closed under $\pi$, Proposition~\ref{HtA} implies that
\[
\lim_{t\rightarrow \infty}  P_{x}^{\pi}(X_t\in A) = P_{x}^{\pi}(\tau_A<\infty ).
\] 
Therefore, 
\[
\lim_{N\rightarrow\infty}\frac{1}{N}\sum_{t=0}^{N-1}P_x^{\pi}(X_t\in A) = P_{x}^{\pi}(\tau_A<\infty ).
\]
\end{pf}

\section{Domain of attraction}\label{domain}

Given a Markov control model $(\mathbb{X},\mathbb{U}, \{\mathbb{U}(x) | x\in\mathbb{X}\}, Q)$, the first problem that we consider is the characterization of the domain of attraction and the escape set of a set $A\subset\SX$, which are defined in Definition \ref{DomAtt&EscSet}. To ensure some stability we assume the following:
\begin{assum}\label{assumA}
Set $A$ is closed under some policy $\pi\in\Pi$. 
\end{assum}
The idea is to characterize both sets in terms of value functions. The first question we ask is whether we can restrict our attention to policies that make $A$ a closed set.
\begin{definition}
Let $\Pi_A\subset \Pi$ be the set of control polices that make $A$ a closed set.
\end{definition}
The following proposition, proved in Appendix \ref{appD}, is fundamental to answer the question above.
\begin{proposition}\label{policyBetter}
Given $\pi\in\Pi$ there exists a policy $\widehat{\pi}\in \Pi_A$ such that for any $x\in\SX$ and $t\geq0$ it holds that $P_x^{\pi}(X_t\in A)\leq P_x^{\widehat{\pi}}(\tau_A\leq t).$
Furthermore, for any $x\in\SX$
\begin{equation*}
\liminf_{t\rightarrow \infty}  P_x^{\pi}(X_t\in A)\leq P_x^{\widehat{\pi}}(\tau_A<\infty ).
\end{equation*}
\end{proposition}
As a corollary, we obtain the following description of $\Lambda_A$ and $\Gamma_A$.
\begin{corollary}\label{enough A closed}
\begin{align*}
\Lambda_A&=\{ x\in \mathbb{X}\mbox{ } | \mbox{ } P_x^{\pi}(\tau_A<\infty )>0 \mbox{ for some policy }\pi\in \Pi_A \}\\ \Gamma_A&=\{ x\in \mathbb{X}\mbox{ } | \mbox{ } P_x^{\pi}(\tau_A<\infty )=0 \mbox{ for all policies }\pi\in \Pi_A \}.
\end{align*}
\end{corollary}
\begin{pf} Clearly, by Proposition \ref{HtA} $\{ x\in \mathbb{X}\mbox{ } | \mbox{ } P_x^{\pi}(\tau_A<\infty )>0 \mbox{ for some policy }\pi\in \Pi_A \} \subset\Lambda_A.$ Let $x\in\Lambda_A$ so there exists a policy $\pi\in\Pi$ such that $\liminf_{t\rightarrow\infty } P_x^{\pi}(X_t\in A )>0$. Let $\widehat{\pi}\in\Pi_A$ be the policy given by the previous proposition. Therefore, $0<\liminf_{t\rightarrow\infty } P_x^{\pi}(X_t\in A )\leq P_x^{\widehat{\pi}}(\tau_A<\infty ).$ 
\end{pf}
Similarly, we can describe $p$-domains of attraction in Definition \ref{pDomAtt} as
\begin{align*}
\Lambda_{A,p}=\left\lbrace x\in \mathbb{X}\mbox{ } \Big| P_x^{\pi}(\tau_A<\infty )\geq p \mbox{ for some policy }\pi\in\Pi_A \right\rbrace.
\end{align*}
Then, we can focus on policies that make $A$ a closed set. In fact, if we consider the average reward function $v^{\pi}$ defined in \eqref{rewardA} and define the value function 
$$V^*(x) =\sup_{\pi\in\Pi_A} v^{\pi}(x),$$ 
then we can use this function to characterize the sets above.
\begin{remark}
Note that $V^*(x)= \sup_{\pi\in \Pi} v^{\pi}(x)$ since any policy $\pi\in\Pi$ can be majorized by the policy $\pi'\in \Pi_A$ defined in Lemma~\ref{NewPolicy}. Indeed, if $x\in A$ it follows clearly and if $x\notin A$, by Proposition \ref{policyBetter} and Lemma \ref{NewPolicy} we obtain that $v^{\pi}(x)\leq v^{\pi'}(x)$. 
\end{remark}
\begin{corollary}\label{DoAvf}
The following characterization of the sets hold:
\begin{align*}
\Lambda_A=\left\lbrace x\in \mathbb{X}\mbox{ } \Big|  V^*(x)>0 \right\rbrace,\\
\Gamma_A =\left\lbrace x\in \mathbb{X}\mbox{ } \Big| V^*(x)=0 \right\rbrace
\end{align*}
and
$$\Lambda_{A,p}=\left\lbrace x\in \mathbb{X}\mbox{ } \Big|  V^*(x)\geq p \right\rbrace.$$
\end{corollary}
\begin{pf}
Let $x\in\Lambda_A$, then Corollary \ref{enough A closed} implies that $P_x^{\pi}(\tau_A<\infty)>0$ for some policy $\pi\in \Pi_A$. Since $A$ is closed under $\pi$, Theorem \ref{eqAv-Ht} implies that $v^{\pi}(x)>0$, as a consequence $v^*(x)>0$. On the other hand, let $x\in\Gamma_A$, so by Corollary \ref{enough A closed} $P_x^{\pi}(\tau_A<\infty)=0$ for all polices $\pi\in\Pi_A$, and the result follows again by Theorem \ref{eqAv-Ht}. Similarly, the result holds for the $p-$domain of attraction.
\end{pf}

\section{Reach-avoid problem}\label{max}

The problem described in the previous section can be seen as a feasibility problem. In this section we will solve a related maximization problem, namely the problem \eqref{P1}. Consider a Markov control model $(\mathbb{X},\mathbb{U}, \{\mathbb{U}(x) | x\in\mathbb{X}\}, Q)$. Let $A,B\subset\mathbb{X}$ be disjoint sets, with their respective hitting times $\tau_A,\tau_{B}$, and an initial distribution $\nu$ over the state space. We will show that this problem is equivalent to a long-run average reward problem with a particular reward function. Our first step will be to rewrite the objective function of \eqref{P1}. To achieve this, let us define a modified stochastic kernel that make sets $A$ and $B$ closed under any policy $\pi$. Given a pair $(x,u)\in\mathbb{K}$ construct the stochastic kernel $\widetilde{Q}$ as follows:
\begin{equation}\label{ModK1}
\begin{cases}
  \widetilde{Q}(A|x,u)=1, \hspace{1cm} \mbox{ if } x\in A\\
  \widetilde{Q}(B|x,u)=1, \hspace{1cm} \mbox{ if } x\in B\\
  \widetilde{Q}(\cdot|x,u)=Q(\cdot|x,u) \hspace{0.5cm} \mbox{ otherwise } \\
\end{cases}
\end{equation}
The measure induced by $\widetilde{Q}$ will be denoted as $\widetilde{P}$ and $\widetilde{\mathbb{E}}$ will denote the expectation with respect to $\widetilde{P}$. By Remark \ref{ClosedQ} the sets $A$ and $B$ are closed for any control policy in $\Pi$. Then, we have the following result proved in Appendix \ref{appB}.

\begin{proposition}\label{propequi}
Given a policy $\pi\in\Pi$, we have that
\begin{equation}\label{equivalent}
P^{\pi}_{\nu}(\tau_A<\tau_B, \tau_A<\infty ) = \widetilde{P}^{\pi}_{\nu}(\tau_A<\infty ).
\end{equation} 
\end{proposition}

Now, we consider the  Markov control model $(\mathbb{X},\mathbb{U}, \{\mathbb{U}(x) | x\in\mathbb{X}\}, \widetilde{Q}, \mathbf{1}_{A})$, that is, the given Markov model with the modified kernel and with the characteristic function of set $A$ as reward function. Let $\widetilde{v}^{\pi}(x)$ be as in \eqref{rewardA} with the modified stochastic kernel.

\begin{thm}\label{ThMarkov}
\begin{align*}
  \sup_{\pi\in\Pi} P^{\pi}_{\nu}(\tau_A<\tau_B,\tau_A<\infty) &= \sup_{\pi\in\Pi_M} \sum_{x\in\mathbb{X}} \widetilde{v}^{\pi}(x)\nu(x)\\
&=\sum_{x\in\mathbb{X}}\widetilde{V}(x)\nu(x),
\end{align*}
where $\widetilde{V}(x)=\sup_{\pi\in\Pi_M}\widetilde{v}^{\pi}(x)$.
\end{thm}
\begin{pf}
Theorem \ref{eqAv-Ht} implies that for $x\in\SX$ and $\pi\in\Pi$, $\widetilde{P}^{\pi}_{\nu}(\tau_A<\infty)=\widetilde{v}^{\pi}(x)$. Hence the result follows from  \eqref{equivalent} and \eqref{eqmarginals}.
\end{pf}

The importance of the result above is that the problem of maximizing the probability of reaching some set $A$ while avoiding a set $B$ can be cast as a long-run average reward problem over stationary Markovian policies.

\section{Reach with hitting constraint}\label{maxcon}

In this section we will solve a constrained version of the previous problem, that is Problem \eqref{P2}. So, again consider a given Markov control model $(\mathbb{X},\mathbb{U}, \{\mathbb{U}(x) | x\in\mathbb{X}\}, Q)$, along with $A,B\subset\mathbb{X}$ disjoint sets and an initial distribution $\nu$. Our objective will be to find a control policy that maximizes the probability of reaching $A$ in such a way that the probability of reaching $B$ is less than some $\epsilon >0$. In this case, however, a Markovian policy might not be the best to solve this problem (recall Theorem \ref{ThMarkov}). The following example shows this fact. To avoid cumbersome notation, in the sequel we will use $\pi$ to denote the stochastic kernels associated with such policy. 

\begin{figure}[b]
\centering
\begin{subfigure}[b]{0.2\textwidth}
\begin{tikzpicture}[->,>=stealth',shorten >=1pt,auto,node distance=1.5cm,
                    thick,main node/.style={circle,draw,font=\sffamily\small\bfseries}]
  \node[main node] (1) {1};
  \node[main node] (2) [right of=1] {2};
  \node[main node] (3) [above of=1] {3};
  \node[main node] (4) [above of=2] {4};
  \node[main node] (5) [above left  of=3] {5};
  \path[every node/.style={font=\sffamily\small}]
    (1) edge [bend left] node {1} (3)
    (2) edge [loop left] node {1} (2)
    (3) edge [bend left] node [left] {0.9} (5)
        edge [bend left] node {0.1} (4)         
    (4) edge [loop above] node {1} (4)     
    (5) edge [loop right] node {1} (5);  
\end{tikzpicture} 
\end{subfigure}
\begin{subfigure}[b]{0.2\textwidth}
\begin{tikzpicture}[->,>=stealth',shorten >=1pt,auto,node distance=1.5cm,
                    thick,main node/.style={circle,draw,font=\sffamily\small\bfseries}]
  \node[main node] (1) {1};
  \node[main node] (2) [right of=1] {2};
  \node[main node] (3) [above of=1] {3};
  \node[main node] (4) [above of=2] {4};
  \node[main node] (5) [above left of=3] {5};
  \path[every node/.style={font=\sffamily\small}]
    (1) edge [bend right] node [right] {1} (3)
    (2) edge [bend right] node {1} (4)
    (3) edge [bend left] node [left] {1} (2)        
    (4) edge [loop above] node {1} (4)    
    (5) edge [loop right] node {1} (5); 
\end{tikzpicture}
\end{subfigure} 
\caption{Control Matrices $u_1,u_2$}
\label{Graph1}
\end{figure}
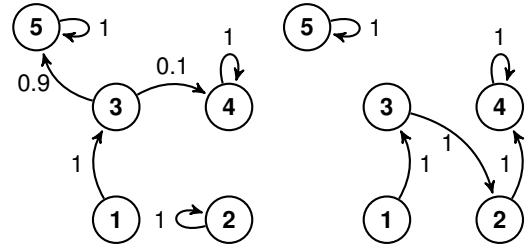

\textbf{Example.} Assume $\mathbb{X}=\{1,2,3,4,5\}$, $\mathbb{U}=\{u_1,u_2\}$ and $\mathbb{U}(x)=\mathbb{U}$ for all $x$. The corresponding control matrices are described in Figure \ref{Graph1}. Consider the sets $A=\{4\},B=\{1,2\}$. Assuming the uniform initial distribution let us consider the problem
\begin{align*}
 & \sup_{\pi\in\Pi}  P^{\pi}_{\nu}(\tau_A<\infty) \\
  &\mbox{s.t. } P^{\pi}_{\nu}(\tau_B<\infty)\leq 0.5
\end{align*}
First of all, regardless of the kind of policy (Markovian or not) we must have that for any policy $\pi$, $P^{\pi}_{1}(\tau_B<\infty)=P^{\pi}_{2}(\tau_B<\infty)=1, P^{\pi}_{4}(\tau_B<\infty)= P^{\pi}_{5}(\tau_B<\infty) = 0$. Thus, to satisfy the restriction we must have that $P^{\pi}_{3}(\tau_B<\infty) \leq 0.5$. Moreover, since $\pi(u_2|3)=P^{\pi}(X_1=2|X_0=3)=P^{\pi}_{3}(\tau_B<\infty)$ then we need that $\pi(u_2|3)\leq 0.5$. We also have that for any policy,
\begin{align*}
   P^{\pi}_{2}(\tau_A<\infty) & =   P^{\pi}(X_1=4|X_0=2)  = \pi(u_2|2), \\
   P^{\pi}_{4}(\tau_A<\infty) & =  1, \quad   P^{\pi}_{5}(\tau_A<\infty) =  0.
\end{align*}
Thus a first obvious choice to maximize $ P^{\pi}_{\nu}(\tau_A<\infty)$ is to select  $\pi(u_2|2)=1$, so that $P^{\pi}_{2}(\tau_A<\infty)=1$. Let $\pi$ be a Markovian policy, so that $P^{\pi}_{1}(\tau_A<\infty) =  P^{\pi}_{3}(\tau_A<\infty)$. If we select $\pi(u_2|3)=1$ we would obtain $P^{\pi}_{1}(\tau_A<\infty) =  P^{\pi}_{3}(\tau_A<\infty)=1$. However, since $\pi(u_2|3)> 0.5$ we can't select such a policy. Moreover, note that for any policy for which $\pi(u_2|3)<1$ we would have $P^{\pi}_{1}(\tau_A<\infty)=P^{\pi}_{3}(\tau_A<\infty)<1$.

On the other hand, if we allow non-Markovian policies we can define a policy for which the value of $P^{\pi}_{3}(\tau_A<\infty)$ is the same as if we have used a Markovian policy, but $P^{\pi}_{1}(\tau_A<\infty)=1$. Indeed, note that 
\begin{align*}
   P&^{\pi}_{1}(\tau_A<\infty)     \\
&= P^{\pi}_{1}(X_1=3,X_2=4) + P^{\pi}_1(X_1=3,X_2=2,X_3=4)\\
   & = \sum_{a_0,a_1\in\mathbb{U}}  P^{\pi}_{1}(a_0,3, a_1,4) 
   +  \sum_{a_0,a_1,a_2\in\mathbb{U}} P^{\pi}_{1}(a_0,3, a_1,2,a_2,4). 
\end{align*}
By conditioning, the terms in the first sum can be written as 
$$Q(4|3,a_1)\pi(a_1|h_1)Q(3|1,a_0)\pi(a_0|1).$$
Similarly, the terms in the second sum can be written as
$$Q(4|2,a_2)\pi(a_2|h_2) Q(2|3,a_1)\pi(a_1|h_1)Q(3|1,a_0)\pi(a_0|1).$$
If we define $\pi(u_2|1)=\pi(u_2|1,u_2,3)=\pi(u_2|1,u_2,3,u_2,2)=1$, we obtain that $ P^{\pi}_{1}(\tau_A<\infty)=1$. Now, we select $\pi(u_2|3)\leq 0.5$ so that the restriction is satisfied. Repeating the same procedure as before we obtain
\begin{align*}
   P&^{\pi}_{3}(\tau_A<\infty)   =   P^{\pi}_{3}(X_1=4) + P^{\pi}_3(X_1=2,X_2=4)\\
   & = \sum_{a_0\in\mathbb{U}}  Q(4|3,a_0)\pi(a_0|3) \\ 
   & + \sum_{a_0,a_1\in\mathbb{U}}  Q(4|2,a_1)\pi(a_1|h_1) Q(2|3,a_0)\pi(a_0|3).
\end{align*}
Thus, by defining $\pi(a_1|h_1)$ in such a way that only depends of the previous state, we can obtain the same value for $ P^{\pi}_{3}(\tau_A<\infty)$ as if we have worked with Markovian policies. Therefore $ P^{\pi}_{\nu}(\tau_A<\infty)$ for such policy is bigger than for any stationary Markovian policy. 

The key point is to realize that if the process $\{X_t\}$ has already gone through the set $B$, then it does not matter if the process hits the set again. But instead, if the chain has not gone through $B$ it is better not to reach the set in order to satisfy the constraint. As a consequence, we will set the problem using control policies that remember whether the process has reached $B$ or not. 

Optimization problems with constraints can be rewritten using Lagrange multiplier. Hence the problem above is equivalent to the following problem 
\[
  \sup_{\pi\in\Pi}  \inf_{\lambda\geq0} \mathcal{L}(\pi,\lambda),
\]
where 
$$\mathcal{L}(\pi,\lambda)= P^{\pi}_{\nu}(\tau_A<\infty) + \lambda (\epsilon - P^{\pi}_{\nu}(\tau_B<\infty)).$$ 
The idea will be to write the Lagrangian function $\mathcal{L}$ as an average reward function with similar ideas as in the previous section.  In order to do this, we consider the set $\widehat{\mathbb{X}} = \mathbb{X}\times \{0,1\}$. Intuitively a state $(x,0)$ indicates that the process has not reached $B$, while a state $(x,1)$ indicates that the process has already reached $B$. Let $\mathbb{U}(x,i)=\mathbb{U}(x)$ for $i=0,1$ and $\widehat{\mathbb{K}}$ the corresponding set of feasible states and actions. We also define the stochastic kernel $\widehat{Q}$ on $\widehat{\mathbb{X}} $ given $\widehat{\mathbb{K}}$ as follows:
\begin{equation}\label{ModK2}
\begin{cases}
  \widehat{Q}((y,0)|(x,0),u)= Q(y|x,u), & \mbox{ if } y\notin B\\
  \widehat{Q}((y,0)|(x,0),u)= 0, &\mbox{ if }  y\in B\\
  \widehat{Q}((y,1)|(x,0),u)= 0, &\mbox{ if } y\notin B\\
  \widehat{Q}((y,1)|(x,0),u)= Q(y|x,u), &\mbox{ if }  y\in B\\
  \widehat{Q}((y,0)|(x,1),u)= 0, &\\
  \widehat{Q}((y,1)|(x,1),u)= Q(y|x,u). &\\
\end{cases}
\end{equation}
Therefore, we have the augmented Markov control model $(\widehat{\mathbb{X}},\mathbb{U}, \{\mathbb{U}(x,i) | (x,i)\in\widehat{\mathbb{X}}\},\widehat{\mathbb{K}},\widehat{Q})$. Let $\widehat{\Pi}$ be the set of policies over the augmented model, $\widehat{\Pi}_M$ the set of Markovian policies, and for any policy $\widehat{\pi}\in\widehat{\Pi}$ we denote by $\widehat{P}^{\widehat{\pi}}$ the measure induced by the policy. The corresponding $\widehat{\mathbb{X}}\times\mathbb{U}$-valued stochastic process is denoted by $\{(X_t,I_t,U_t)\}_{t\geq0}$. Histories in this model are denoted by $\widehat{h}_t=(h_t,\overline{i}^t)$, with $\overline{i}^t$ the history up to time $t$ of the process $\{I_t\}$.

\begin{definition}\label{indpi}
Given a $\pi\in\Pi$ we define a policy $\widehat{\pi}\in\widehat{\Pi}$ as follows: For any $t\geq0$
$$\widehat{\pi}(u| \widehat{h}_t)=\pi(u|h_t),$$
whenever $\widehat{h}_t$ satisfies that $\overline{x}^t\notin B$ and $ \overline{i}^t\notin \{1\}$, or there is some $0\leq s\leq t$ such that $\overline{x}^{s-1}\notin B$, $x_s\in B$, $\overline{i}^{s-1}\notin\{1\}$ and $i_r=1$ for $s\leq r\leq t$. Otherwise, we define $\widehat{\pi}(u| \widehat{h}_t)=\delta_{u_{0}}$ for any $u_0\in\mathbb{U}(x_t)$. Note that $s=0$ implies that $\overline{i}^{t}\notin\{0\}$.
\end{definition}

The next result allows to express joint distributions of the original model in terms of joint distributions of the augmented model. Its proof can be found in Appendix \ref{appLemma}.
\begin{lemma}\label{LemmaRestr}
Let $\pi\in\Pi$ and consider the policy $\widehat{\pi}\in\widehat{\Pi}$ of the previous definition. Given $h_t\in H_t$ and $u_t\in\mathbb{U}(x_t)$, let $\widehat{h}_t =((x, 0) , u_0 , \ldots , (x_{t-1},0), u_{t-1}, (x_{t},0) )$ if $\overline{x}^t\notin B$, and $\widehat{h}_t=((x, 0), \ldots , (x_{s-1},0), u_{s-1}, (x_{s},1),\ldots , (x_{t},1) )$ if $\overline{x}^{s-1}\notin B,x_s \in B$. Then, 
\[
P^{\pi}_x(h_t, u_t) =
\begin{cases}
 \widehat{P}^{\widehat{\pi}}_{(x,0)} (\widehat{h_t} , u_t ) \mbox{ if } x\notin B \\  
  \widehat{P}^{\widehat{\pi}}_{(x,1)} (\widehat{h_t} , u_t ) \mbox{ if } x\in B.
\end{cases}
\]
\end{lemma}

Note that by the definition of $\widehat{Q}$, the set $\SX\times\{1\}$ is closed under every policy in the augmented Markov model. Also, as in Section \ref{max}, we can redefine the kernel $\widehat{Q}$ to make the set $A\times\{0,1\}$ closed under any policy, that is, $\widehat{Q}(A\times\{i\}|(x,i),u)= 1$ if $x\in A$ and $i=0,1$. With this redefined kernel we consider the Markov model along with the reward function
$$r(x,i)=(\mathbf{1}_{A\times \{0,1\}} - \lambda\mathbf{1}_{\SX\times \{1\}})(x,i).$$
Given a policy $\widehat{\pi}\in\widehat{\Pi}$ and $(x,i)\in\widehat{\mathbb{X}}$, we consider the average reward function given by
\begin{equation} \label{avcost}
 v^{\widehat{\pi}}(x,i) :=\limsup_{N\rightarrow\infty}\frac{1}{N}\widehat{\mathbb{E}}_{(x,i)}^{\widehat{\pi}}\Big[ \sum_{t=0}^{N-1}[\mathbf{1}_{A\times \{0,1\}}  - \lambda\mathbf{1}_{\SX\times \{1\}}](X_t,I_t) \Big]. 
\end{equation}

The next theorem proved in Appendix \ref{appC} shows that the Lagrangean $\mathcal{L}$ can be written in terms of this function.
\begin{thm}\label{Linaverage}
Let $\pi\in\Pi$ and $\widehat{\pi}\in\widehat{\Pi}$ the policy defined in Definition \ref{indpi}. Let $x\in\SX$, then
\begin{align*}
P^{\pi}_{x}(\tau_A<\infty) + \lambda &\left(\epsilon - P^{\pi}_{x}(\tau_B<\infty)\right)= \\
 &
\begin{cases}
  v^{\widehat{\pi}}(x,0) + \lambda\epsilon,  \mbox{ if } x\notin B \\
 v^{\widehat{\pi}}(x,1) + \lambda\epsilon,  \mbox{ if } x\in B.
\end{cases}\\
\end{align*}
\end{thm}
Now, since Definition \ref{indpi} does not necessarily recover all policies in $\widehat{\Pi}$, the previous theorem shows that the following is an upper bound for \eqref{P2}
\begin{equation}\label{upper}
\sup_{\widehat{\pi}\in\widehat{\Pi}}  \inf_{\lambda\geq0} \lambda\epsilon+ \sum_{x\notin B} v^{\widehat{\pi}}(x,0)\nu(x)+\sum_{x\in B}v^{\widehat{\pi}}(x,1)\nu(x).
\end{equation}
Also, note that in the problem above we can consider only Markovian policies by \eqref{eqmarginals}. Furthermore, given a policy $\widehat{\pi}\in\widehat{\Pi}_M$ and stationary we can define a policy $\pi\in\Pi$ by
\begin{equation}\label{aumentado_a_original}
 \pi(u|h_t) =\begin{cases}
\widehat{\pi}(u|x_t,0) & \mbox{ if } \overline{x}^{t}\notin B \\
\widehat{\pi}(u|x_t,1) & \mbox{ otherwise. }
\end{cases}
\end{equation}
Let the set of policies obtained by \eqref{aumentado_a_original} be denoted as $\Pi_B\subset\Pi$, that is, policies that  do not depend on the whole history but only whether the process has gone through set $B$ or not. Hence, Definition \ref{indpi} and Equation \eqref{aumentado_a_original} define a one-to-one relation between $\Pi_B$ and $\widehat{\Pi}_M$.
\begin{thm}\label{corP2}
Problem \eqref{P2} is equivalent to
\begin{align}\label{eqP2b}
  &\sup\limits_{\pi\in\Pi_B}  P^{\pi}_{\nu}(\tau_A<\infty) \\\notag
  &\mbox{s.t. }  P^{\pi}_{\nu}(\tau_B<\infty)\leq \epsilon,
\end{align}
its optimal value is equal to
\begin{equation}\label{eqP2}
\sup_{\widehat{\pi}\in\widehat{\Pi}_M}  \inf_{\lambda\geq0} \lambda\epsilon+ \sum_{x\notin B} v^{\widehat{\pi}}(x,0)\nu(x)+\sum_{x\in B}v^{\widehat{\pi}}(x,1)\nu(x),
\end{equation}
and any optimal policy of \eqref{eqP2} produces an optimal policy of \eqref{P2} through \eqref{aumentado_a_original}.
\end{thm}
\begin{pf}
Let $P^*$ be the optimal value of Problem \eqref{P2}. Then $P^*$ is bounded above by \eqref{upper} which has the same value as \eqref{eqP2}, which by \eqref{aumentado_a_original} has the same value as the optimal value of \eqref{eqP2b}, which is bounded above by $P^*$. Hence all these expressions are equivalent.
\end{pf}

\section{Linear programming formulations: Finite case}\label{lpFinite}

In this section we consider the case where the state space $\mathbb{X}$ and the action space $\mathbb{U}$ are finite and present linear programs that solve our problems. We know from Theorem \ref{AverageDual} and Corollary \ref{DoAvf} that in order to find the $p$-domain of attraction of a given set $A$, we need to solve the following linear program
\begin{align}\label{ReaP}\tag{ReaP}
  \min & \sum_{x\in\mathbb{X}} v(x) \\\notag
  \mbox{s.t. }& \\\notag
 & v(x) \geq \sum_{j\in\mathbb{X}} P_{x,x'}^u v(x'),\\\notag
 & v(x)  \geq  \mathbf{1}_A(x)+ \sum_{x'\in\mathbb{X}} P_{x,x'}^u h(x') - h(x), \\\notag
&\forall x\in\mathbb{X},u\in\mathbb{U}(x)
\end{align}
where $P_{x,j}^u=Q(j|x,u)$. Similarly, in order to solve Problem \eqref{P1}, from Theorem \ref{ThMarkov} we need to solve the same linear program as above with $\widetilde{P}_{x,j}^u=\widetilde{Q}(j|x,u)$ instead of $P_{x,j}$, where $\widetilde{Q}$ is defined in \eqref{ModK1}.  In both cases, an optimal stationary Markovian policy can be found by solving the corresponding dual problems, as in Theorem \ref{AverageDual}. The dual problem is the following:
\begin{align}\label{ReaDual}\tag{ReaD}
  &\max \sum_{x\in A} \sum_{u\in\mathbb{U}(x)}\alpha(x,u)  \\\notag
  &\mbox{s.t.} \\\notag
 & \sum_{u\in\mathbb{U}(x)} \alpha(x,u) - \sum_{x'\in\mathbb{X}\atop u\in\mathbb{U}(x')}  P_{x',x}^u \alpha(x',u)  = 0, \\\notag
 & \sum_{u\in\mathbb{U}(x)} \alpha(x,u)  +   \beta(x,u) - \sum_{x'\in\mathbb{X}\atop u\in\mathbb{U}(x')} P_{x',x}^u \beta(x',u)  = 1, \\\notag
 & \alpha(x,u)\geq 0 ,\beta(x,u)\geq 0,\quad\forall x\in\mathbb{X},u\in\mathbb{U}(x).
\end{align}

\begin{remark} Note that in both cases the initial distribution $\nu$ does not play any role in order to find either the value functions $V^*$ and $\widetilde{V}$ and the optimal policies.
\end{remark}

For Problem \eqref{P2}, the linear programming formulation is not as straightforward as in the previous problems. In particular, we would like to switch the $\inf$ with the $\sup$ in \eqref{eqP2}, that is, we would like to show that it is equivalent to its dual problem. Recall that \eqref{P2} can be written as
$$P^*=\sup_{\pi\in\Pi_B}  \inf_{\lambda\geq0} P^{\pi}_{\nu}(\tau_A<\infty) + \lambda (\epsilon - P^{\pi}_{\nu}(\tau_B<\infty)),$$
which is bounded above by the optimal value of its dual problem
$$D^*=\inf_{\lambda\geq0} \sup_{\pi\in\Pi_B}  P^{\pi}_{\nu}(\tau_A<\infty) + \lambda (\epsilon - P^{\pi}_{\nu}(\tau_B<\infty)).$$
By Theorem \ref{corP2} we can write the above problem as
\begin{align}\label{D1}\tag{D1}
  &\inf_{\lambda\geq0} \lambda\epsilon+\max \sum_{j\in A\times\{0,1\},\atop u\in\mathbb{U}(j)} \alpha(j,u) -\lambda\sum_{j\in \mathbb{X}\times\{1\}\atop u\in\mathbb{U}(j)} \alpha(j,u) \\\notag
 & \mbox{s.t.}\\\notag
 & \sum_{u\in\mathbb{U}(j)} \alpha(j,u) -  \sum_{j'\in\widehat{\mathbb{X}}\atop u\in\mathbb{U}(j')}  \widehat{P}_{j',j}^u \alpha(j,u)  = 0, \\\notag
 & \sum_{u\in\mathbb{U}(j)} \alpha(j,u)  +   \beta(j,u) - \sum_{j'\in\widehat{\mathbb{X}}\atop u\in\mathbb{U}(j')} \widehat{P}_{j',j}^u \beta(j',u)  = \widehat{\nu}(j), \\\notag
 & \alpha(j,u)\geq 0 ,\beta(j,u)\geq 0,\quad\forall j\in\widehat{\mathbb{X}},u\in\mathbb{U}(j).
\end{align}
 with $\widehat{P}_{(x,i),(x',i')}^u=\widehat{Q}((x',i')|(x,i),u)$,  where $\widehat{Q}$ is defined in \eqref{ModK2}, and 
$$\widehat{\nu}(x,i)=
\begin{cases}
  \nu(x), & \mbox{ if } x\notin B,i=0\\
  \nu(x), &\mbox{ if }  x\in B,i=1\\
  0, &\mbox{ otherwise.} 
\end{cases}
$$
Now, by strong duality of linear programming we further obtain that 
\begin{align}\label{D2}\tag{D2}
 D^*&\\\notag
=&\max \sum_{j\in A\times\{0,1\},\atop u\in\mathbb{U}(j)} \alpha(j,u)  \\\notag
  &\mbox{s.t.} \\\notag
 & \sum_{u\in\mathbb{U}(j)} \alpha(j,u) -  \sum_{j'\in\widehat{\mathbb{X}}\atop u\in\mathbb{U}(j')}  \widehat{P}_{j',j}^u \alpha(j,u)  = 0, \\\notag
 & \sum_{u\in\mathbb{U}(j)} \alpha(j,u)  +   \beta(j,u) - \sum_{j'\in\widehat{\mathbb{X}}\atop u\in\mathbb{U}(j')} \widehat{P}_{j',j}^u \beta(j',u)  = \widehat{\nu}(j), \\\notag
&\sum_{j\in \mathbb{X}\times\{1\}\atop u\in\mathbb{U}(j)} \alpha(j,u)\leq \epsilon,\\\notag
 & \alpha(j,u)\geq 0 ,\beta(j,u)\geq 0,\quad\forall j\in\widehat{\mathbb{X}},u\in\mathbb{U}(j).
\end{align}
If Problem \eqref{D2} is infeasible, then $D^*=-\infty$ and therefore $P^*=-\infty$, that is, Problem \eqref{P2} is infeasible. On the other hand, if \eqref{D2} is finite (note that it cannot be unbounded) with optimal solution $(\alpha^*,\beta^*$), the $\inf$ over $\lambda\geq0$ in \eqref{D1} is attained at some $\lambda^*$ such that
$$\lambda^*\left(\epsilon-\sum_{x\in \mathbb{X},u\in\mathbb{U}(x)} \alpha^*((x,1),u)\right)=0,$$
by Complementary Slackness condition.
Let $\widehat{\pi}^*\in\widehat{\Pi}_M$ the stationary optimal policy induced by $(\alpha^*,\beta^*)$ given by Theorem \ref{AverageDual}, and $\pi^*\in\Pi_B$ its corresponding over $\mathbb{X}$ given by \eqref{aumentado_a_original}. Therefore,
\begin{align*}
P^*&=\sup_{\pi\in\Pi_B}  \inf_{\lambda\geq0} P^{\pi}_{\nu}(\tau_A<\infty) + \lambda (\epsilon - P^{\pi}_{\nu}(\tau_B<\infty))\\
&\geq\inf_{\lambda\geq0} P^{\pi^*}_{\nu}(\tau_A<\infty) + \lambda (\epsilon - P^{\pi^*}_{\nu}(\tau_B<\infty))\\
&=\inf_{\lambda\geq0} \lambda\epsilon +\sum_{j\in\widehat{\SX}}v^{\widehat{\pi}^*}(j)\widehat{\nu}(j)\\
&=\inf_{\lambda\geq0} \lambda\epsilon+\sum_{j\in A\times\{0,1\},\atop u\in\mathbb{U}(j)} \alpha^*(j,u) -\lambda\sum_{j\in \mathbb{X}\times\{1\}\atop u\in\mathbb{U}(j)} \alpha^*(j,u)\\
&=\lambda^*\epsilon+\sum_{j\in A\times\{0,1\},\atop u\in\mathbb{U}(j)} \alpha^*(j,u) -\lambda\sum_{j\in \mathbb{X}\times\{1\}\atop u\in\mathbb{U}(j)} \alpha^*(j,u)\\
&=\sum_{j\in A\times\{0,1\},\atop u\in\mathbb{U}(j)} \alpha^*(j,u)=D^*\geq P^*.
\end{align*}
Hence, we just proved the following result.

\begin{thm}\label{finite_case}
Suppose $\mathbb{X},\mathbb{U}$ are finite. Then \eqref{P2} satisfies strong duality, that is, it is equivalent to 
\begin{align*}
\min_{\lambda \geq 0} \max_{\pi\in\Pi} \mathcal{L}(\pi,\lambda).
\end{align*}
Furthermore, it can be solved by the linear program \eqref{D2} in the augmented model and recover an optimal policy in $\Pi_B$ through  \eqref{aumentado_a_original}.
\end{thm}

\section{Numerical example}\label{example}

\begin{figure}[h!]
\centering
\includegraphics[scale=0.6]{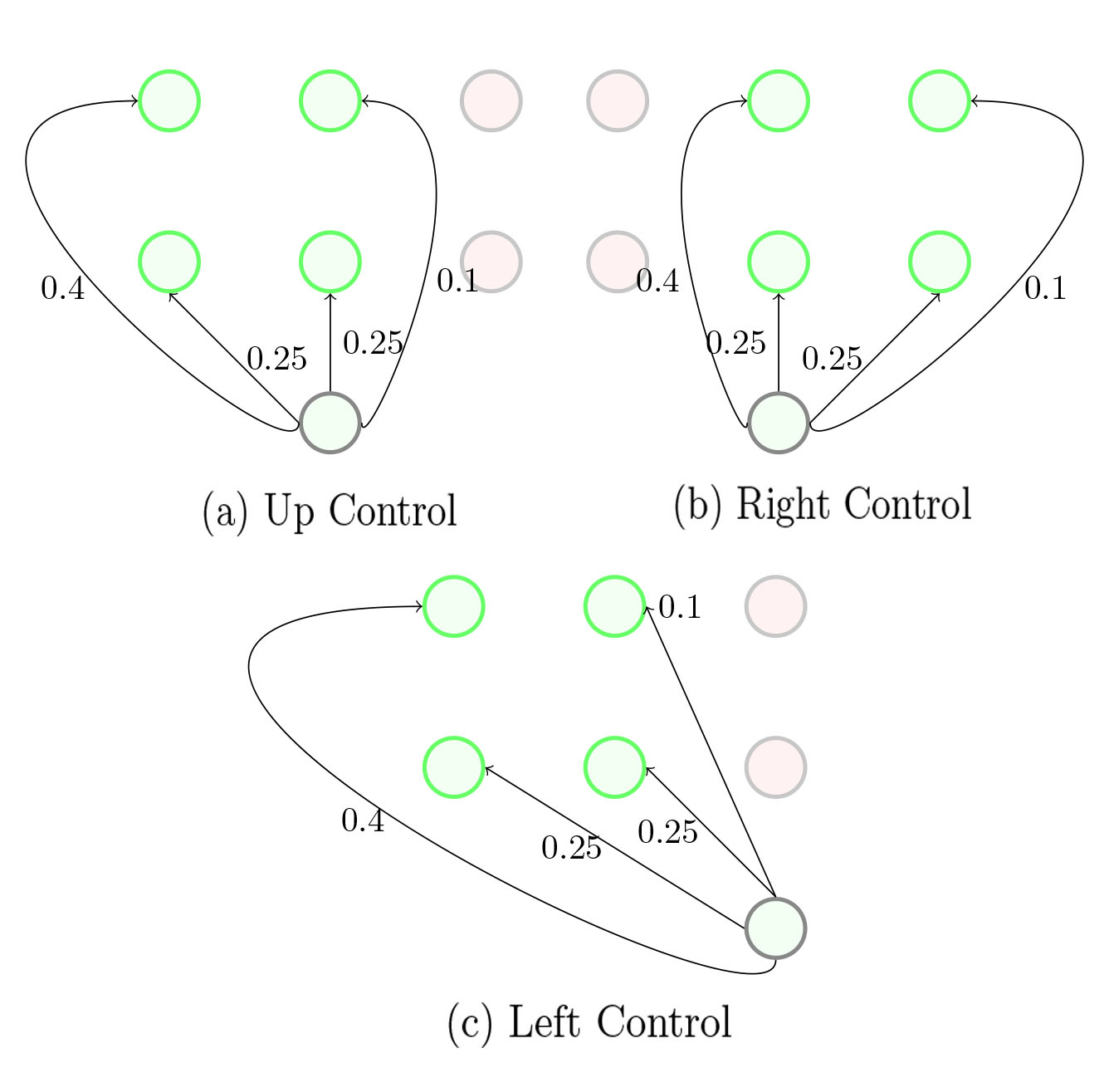}
\caption{Controls}\label{control}
\end{figure}

To illustrate our results we consider an object that navigates over a grid under the influence of a north-west wind. The object has three controls available as shown in Figure \ref{control}. We assume that the states at the upper boundary of the grid are absorbing states and some adjustments are done in the left and right boundaries so that the object does not leave the grid. This example is similar to the Zermelo Navigation problem presented in \cite{Esfahani} in the context of continuous time problems.

\subsection{Domain of attraction}

\begin{figure}[h!]
\centering
   \begin{subfigure}[b]{.5\linewidth}
      \includegraphics[width=1\textwidth]{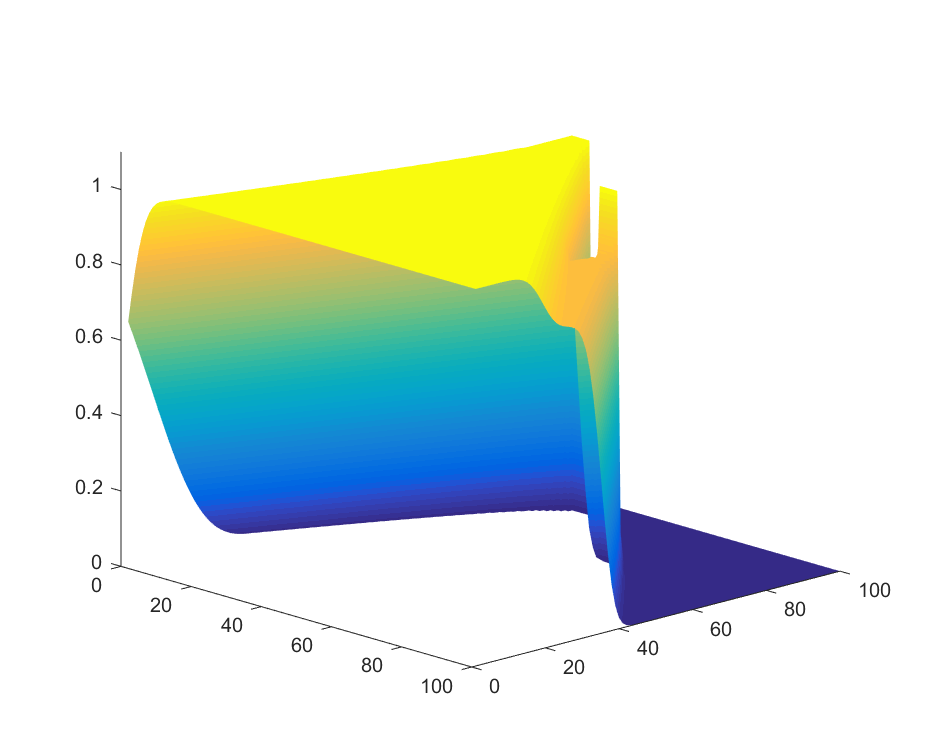}
\caption{Surface of $V^*(x)$.\\}
   \end{subfigure}
\begin{subfigure}[b]{.45\linewidth}
      \includegraphics[width=1\textwidth]{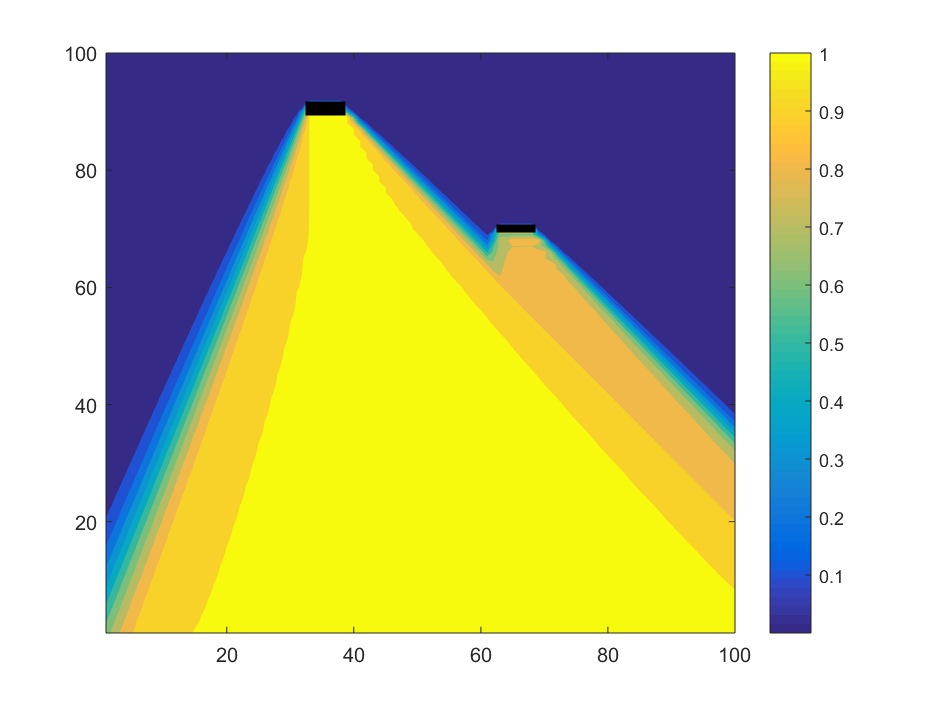}
\caption{Level sets of $V^*(x)$.}
   \end{subfigure}
      \caption{$p$-domains $\Lambda_p$ and escape set $\Gamma$.}\label{p-domains}
\centering
\end{figure}

In order to show the findings of Section \ref{domain}, we consider a 100 by 100 grid with a closed set $A$ in the central region of the grid  marked with black squares. Figure \ref{p-domains} shows the surface and level sets of the function $V^*(x)$ which defines the $p$-domains $\Lambda_p$. The scape set $\Gamma$ corresponds to the states with value function equal to zero. This function was found by solving the linear program \eqref{ReaP}. We note that no state outside of $A$ belongs to $\Lambda_1$. 

\subsection{Reach and avoid}

\begin{figure}[h!]
\centering
\includegraphics[scale=0.33]{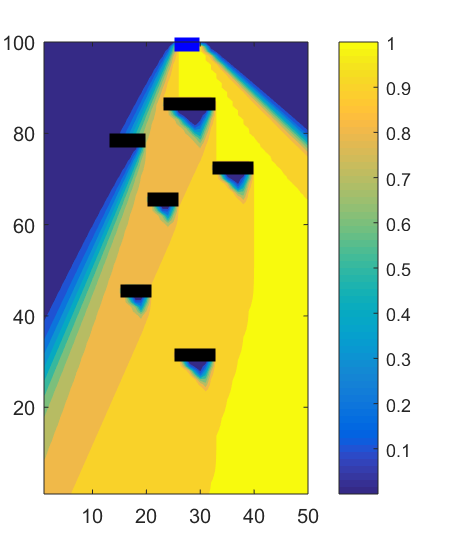}
\caption{Level sets of $\widetilde{V}(x)$.}\label{P1vf}
\end{figure}

\begin{figure}[h!]
\centering
   \begin{subfigure}[b]{.33\linewidth}
      \includegraphics[width=1\textwidth]{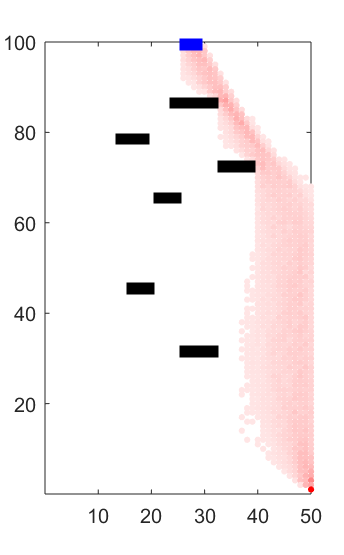}
\caption{Initial state (50,1)}
   \end{subfigure}
\begin{subfigure}[b]{.33\linewidth}
      \includegraphics[width=1\textwidth]{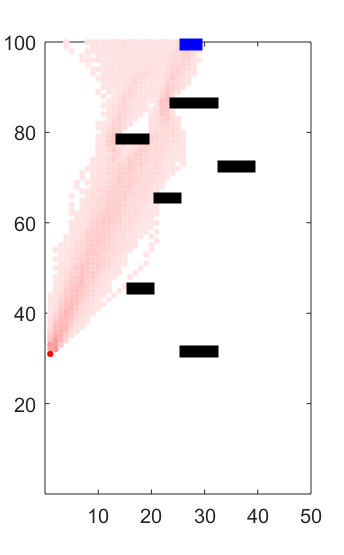}
\caption{Initial state (1,30)}
   \end{subfigure}
      \caption{Trajectories.}\label{P1traj}
\centering
\end{figure}

For Problem \eqref{P1} we consider a 100 by 50 grid with the set $A$ a portion of the upper boundary of the grid, marked with dark blue squares, and set $B$ a number of obstacles spread over the grid, marked with black squares. Figure \ref{P1vf} shows the level sets of the function $\widetilde{V}(x)$ computed by solving the linear program \eqref{ReaP}. In Figure \ref{P1traj} we show the paths of 500 simulated trajectories of the object under the optimal policy obtained from the linear program \eqref{ReaDual}. We choose two different starting states from different level sets according to Figure \ref{P1vf}. Note that inthe first case all trajectories hit the target set $A$, while in the second case most of them drift away from the set. This situation agrees with the Figure \ref{P1vf}.

\subsection{Reach with constraint}

\begin{figure}[h!]
\centering
   \begin{subfigure}[b]{.38\linewidth}
      \includegraphics[width=1\textwidth]{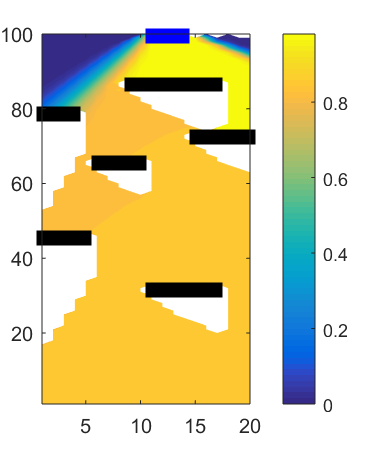}
\caption{$\epsilon=0.01$}
   \end{subfigure}
\begin{subfigure}[b]{.38\linewidth}
      \includegraphics[width=1\textwidth]{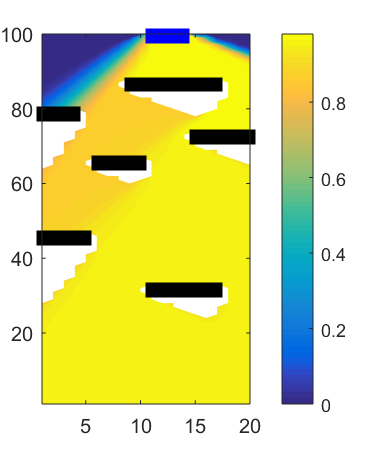}
\caption{$\epsilon=0.2$}
   \end{subfigure}
\begin{subfigure}[b]{.38\linewidth}
      \includegraphics[width=1\textwidth]{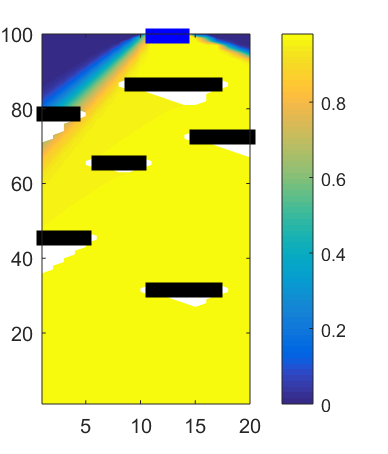}
\caption{$\epsilon=0.8$}
   \end{subfigure}
      \caption{Optimal value of \eqref{P2}.}\label{P2value}
\centering
\end{figure}

\begin{figure}[t!]
\centering
   \begin{subfigure}[b]{.26\linewidth}
      \includegraphics[width=1\textwidth]{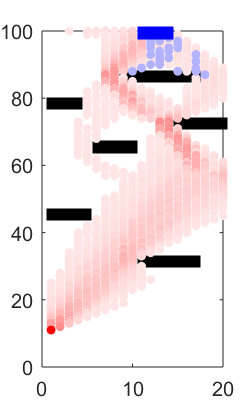}
   \end{subfigure}
\begin{subfigure}[b]{.26\linewidth}
      \includegraphics[width=1\textwidth]{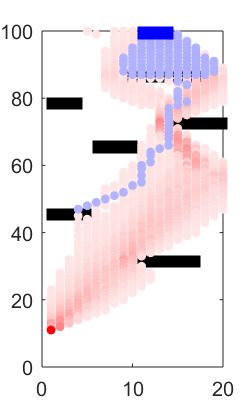}
   \end{subfigure}
\begin{subfigure}[b]{.26\linewidth}
      \includegraphics[width=1\textwidth]{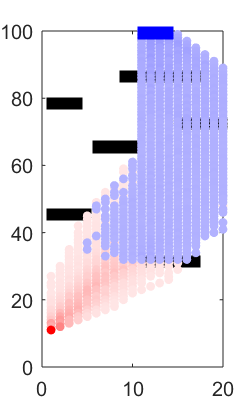}
   \end{subfigure}

   \begin{subfigure}[b]{.26\linewidth}
      \includegraphics[width=1\textwidth]{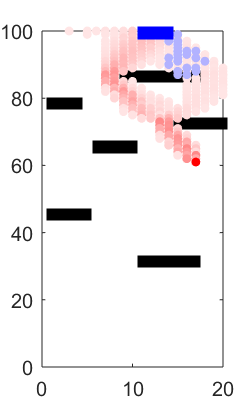}
\caption{$\epsilon=0.01$}
   \end{subfigure}
\begin{subfigure}[b]{.26\linewidth}
      \includegraphics[width=1\textwidth]{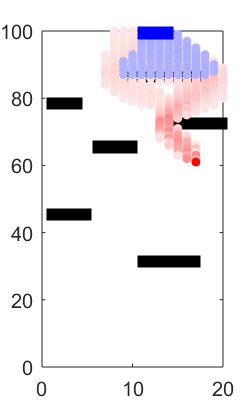}
\caption{$\epsilon=0.2$}
   \end{subfigure}
\begin{subfigure}[b]{.26\linewidth}
      \includegraphics[width=1\textwidth]{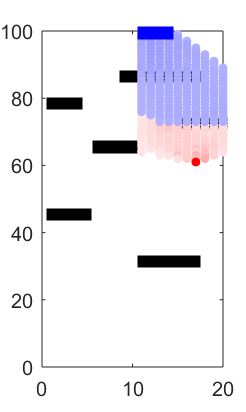}
\caption{$\epsilon=0.8$}
   \end{subfigure}
      \caption{Trajectories.}\label{P2traj}
\centering
\end{figure}

For Problem \eqref{P2} we consider a 100 by 20 grid and sets $A$ and $B$ marked with dark blue and black squares, respectively. It is important to note that the initial distribution $\nu$ plays a key role in the feasibility of the problem. Figure \ref{P2value} shows the level sets of the function
$$\widehat{V}(x):= \max \limits_{\pi\in\Pi_B} P^{\pi}_x(\tau_A<\infty)\quad\text{s.t}\quad P^{\pi}_x(\tau_B<\infty)\leq\epsilon,$$
for different values of $\epsilon$. White regions represents the states for which the problem above is infeasible. As expected, the number of infeasible states decrease with bigger values of $\epsilon$. In Figure \ref{P2traj} we show the paths of 500 trajectories under the optimal policies with different initial distributions and values of $\epsilon$. Red trajectories correspond to optimal trajectories before hitting set $B$ and blue trajetories after hitting this set. Note also that blue trajectories do not avoid the obstacles since they were already hit.

\section{Future work}\label{future}

There are various directions for future research. The first one is to apply the ideas presented in this work to a general context that includes general state and actions spaces and continuous time controlled Markov process. The second one is to study the robust counterpart of these reachability problems. The recent and increasing literature in robust MDPs using different classes of ambiguity sets can be applied in our context. Finally, it would be interesting to to consider problems with moving target and obstacle sets.

\begin{ack}                              
The authors were partially supported by the FAPA funds from Universidad de los Andes.  
\end{ack}

\bibliographystyle{plain}

\bibliography{biblio}

\begin{bibdiv}
\begin{biblist}

\bib{abate}{article}{
      author={Abate, Alessandro},
      author={Prandini, Maria},
      author={Lygeros, John},
      author={Sastry, Shankar},
       title={Probabilistic reachability and safety for controlled discrete
  time stochastic hybrid systems},
        date={2008},
        ISSN={0005-1098},
     journal={Automatica},
      volume={44},
      number={11},
       pages={2724 \ndash  2734},
  url={http://www.sciencedirect.com/science/article/pii/S0005109808002677},
}

\bib{arevalo}{inproceedings}{
      author={{Arvelo}, E.},
      author={{Kim}, E.},
      author={{Martins}, N.~C.},
       title={Maximal persistent surveillance under safety constraints},
        date={2013},
   booktitle={2013 {IEEE} {I}nternational {C}onference on {R}obotics and
  {A}utomation},
       pages={4048\ndash 4053},
}

\bib{arvelo2883control}{article}{
      author={Arvelo, E},
      author={Martins, NC},
       title={Maximizing the set of recurrent states of an mdp, subject to
  convex constraints},
        date={2014},
     journal={Automatica},
      volume={50},
      number={3},
       pages={994\ndash 998},
}

\bib{bertsekas1995dynamic}{book}{
      author={Bertsekas, Dimitri~P},
       title={Dynamic programming and optimal control},
   publisher={Athena Scientific},
        date={1995},
      volume={1},
}

\bib{bertsekas1978stochastic}{book}{
      author={Bertsekas, Dimitri~P},
      author={Shreve, Steven~E},
       title={Stochastic optimal control: The discrete time case},
   publisher={Athena Scientific},
        date={1996},
}

\bib{Bouchard}{article}{
      author={Bouchard, Bruno},
      author={Elie, Romuald},
      author={Touzi, Nizar},
       title={Stochastic target problems with controlled loss},
        date={2010},
     journal={SIAM Journal on Control and Optimization},
      volume={48},
      number={5},
       pages={3123\ndash 3150},
}

\bib{camilli2008control}{article}{
      author={Camilli, Fabio},
      author={Gr{\"u}ne, Lars},
      author={Wirth, Fabian},
       title={Control {L}yapunov functions and {Z}ubov's method},
        date={2008},
     journal={SIAM Journal on Control and Optimization},
      volume={47},
      number={1},
       pages={301\ndash 326},
}

\bib{camilli2006zubov}{article}{
      author={Camilli, Fabio},
      author={Loreti, Paola},
       title={A {Z}ubov's method for stochastic differential equations},
        date={2006},
     journal={Nonlinear Differential Equations and Applications NoDEA},
      volume={13},
      number={2},
       pages={205\ndash 222},
}

\bib{chatterjee2011maximizing}{article}{
      author={Chatterjee, Debasish},
      author={Cinquemani, Eugenio},
      author={Lygeros, John},
       title={Maximizing the probability of attaining a target prior to
  extinction},
        date={2011},
     journal={Nonlinear Analysis: Hybrid Systems},
      volume={5},
      number={2},
       pages={367\ndash 381},
}

\bib{Esfahani}{article}{
      author={Esfahani, Peyman~Mohajerin},
      author={Chatterjee, Debasish},
      author={Lygeros, John},
       title={The stochastic reach-avoid problem and set characterization for
  diffusions},
        date={2016},
        ISSN={0005-1098},
     journal={Automatica},
      volume={70},
       pages={43 \ndash  56},
}

\bib{hahn1967stability}{book}{
      author={Hahn, Wolfgang},
      author={Baartz, Arne~P},
       title={Stability of motion},
   publisher={Springer},
        date={1967},
      volume={138},
}

\bib{hernandez2012further}{book}{
      author={Hern{\'a}ndez-Lerma, On{\'e}simo},
      author={Lasserre, Jean~B},
       title={Further topics on discrete-time markov control processes},
   publisher={Springer Science \& Business Media},
        date={2012},
      volume={42},
}

\bib{hernandez1996discrete}{book}{
      author={Hern{\'a}ndez-Lerma, On{\'e}simo},
      author={Lasserre, Jean~Bernard},
       title={Discrete-time markov control processes: basic optimality
  criteria},
   publisher={Springer},
        date={1996},
}

\bib{puterman2014markov}{book}{
      author={Puterman, Martin~L},
       title={Markov decision processes: discrete stochastic dynamic
  programming},
   publisher={John Wiley \& Sons},
        date={2014},
}

\bib{Soner}{article}{
      author={Soner, H.~Mete},
      author={Touzi, Nizar},
       title={Stochastic target problems, dynamic programming, and viscosity
  solutions},
        date={2002},
     journal={SIAM Journal on Control and Optimization},
      volume={41},
      number={2},
       pages={404\ndash 424},
}

\bib{summers2010verification}{article}{
      author={Summers, Sean},
      author={Lygeros, John},
       title={Verification of discrete time stochastic hybrid systems: A
  stochastic reach-avoid decision problem},
        date={2010},
     journal={Automatica},
      volume={46},
      number={12},
       pages={1951\ndash 1961},
}

\bib{zubov1964methods}{book}{
      author={Zubov, Vladimir~Ivanovich},
       title={Methods of {AM} {L}yapunov and their application},
   publisher={P. Noordhoff},
        date={1964},
}

\end{biblist}
\end{bibdiv}

\appendix
\section{Proof of Proposition \ref{HtA}}\label{appA}
To prove the proposition we will first need the two following lemmas.
\begin{lemma}\label{nullevent} 
Let $A$ be closed for a policy $\pi\in\Pi$ and consider the event $E=\{ \exists s,t \in \mathbb{N} \mbox{ such that } s\leq t, X_s\in A, X_t\notin A\}$. Then $P_x^{\pi}(E)=0$ for any $x\in \mathbb{X}$.
\end{lemma}
\begin{pf} We can write such an event as $E=\bigcup_{s}\bigcup_{s < t }\{ X_s\in A, X_t\notin A\}$. Therefore, $P_x^{\pi}(E) \leq \sum_s \sum_{s < t} P^{\pi}_x( X_s\in A , X_t \notin A).$ Let us see that for all $s < t$, $P_x^{\pi}(X_s\in A , X_t \notin A)=0$. We have that, 
\begin{align*}
   P_x^{\pi}&(X_s\in A , X_t \notin A) \\
& =  \sum_{x_i\in\mathbb{X},x_s\in A, u_i \in\mathbb{U}(x_i)}  P_x^{\pi}(u_0, x_1, u_1,\ldots,x_{t-1}, X_t \notin A )\\
 & =  \sum_{ x_i\in \mathbb{X},x_s\in A, u_i \in\mathbb{U}(x_i)} P^{\pi}(X_t \notin A  | h_{t-1} )P^{\pi}(h_{t-1}|x) =  0,
\end{align*}
where last equality follows from the fact that $A$ is closed under $\pi$, and so $P^{\pi}(X_t \notin A  | h_{t-1} )=0$.
\end{pf} 

\begin{lemma}\label{Ht<n=X_n in A}
Let $A$ be closed for a policy $\pi\in\Pi$. Then, for any $x\in \mathbb{X}$
\[
P_x^{\pi}(\tau_A\leq t) = P_x^{\pi}(X_t\in A).
\] 
\end{lemma}

\begin{pf} Define $B:=\{ \omega \in \Omega \hspace{0.2cm} | \hspace{0.2cm} \forall \hspace{0.2cm} m > \tau_A(\omega) , X_m(\omega)\in A\}$. So that, 
\[
\{\tau_A\leq t\} = (\{\tau_A\leq t\}\cap B ) \cup (\{\tau_A\leq t\}\cap B^c)
\]
and
\[
\{X_t\in A\} = (\{X_t\in A\}\cap B ) \cup (\{X_t\in A\}\cap B^c).
\]
Let us see that the events $\{X_t\in A\}\cap B $ and $\{\tau_A\leq t\}\cap B$ are equal. It is clear that $\{X_t\in A\}\cap B\subset\{\tau_A\leq t\}\cap B$. Now, let $\omega\in \{\tau_A\leq t\}\cap B$ so $\tau_A(\omega)\leq t$. Since $\omega\in B$, for all $m\geq \tau_A(\omega)$ we have $X_m(\omega)\in A$, in particular $X_t(\omega)\in A$, so that $\{\tau_A\leq t\}\cap B\subset \{X_t\in A\}\cap B$. Therefore,
\[P_x^{\pi}( \{\tau_A\leq t\}\cap B) = P_x^{\pi}( \{X_t\in A\}\cap B).\]
Finally, consider the event $E=\{ \exists s,t\in \mathbb{N}\mbox{ such that } j\leq m, X_j\in A, X_m\notin A \}$. Note that, $\{\tau_A\leq t\}\cap B^c\subset E$ and $\{X_t\in A\}\cap B^c\subset E$. Since $A$ is closed, Lemma~\ref{nullevent} implies
$P_x(E)^{\pi}=0$, so that,
\[
 0=P_x^{\pi}(\{\tau_A\leq t \} \cap B^c)=P_x^{\pi}(\{X_t\in A\} \cap B^c ).
\]
Combining both results we obtain the required equality. 
\end{pf}
Now, the proposition follows from the previous lemma and Monotone Convergence Theorem.

\section{Proof of Proposition \ref{policyBetter}}\label{appD}

In order to prove the proposition we first need the following lemma.

\begin{lemma}\label{NewPolicy}
Given $\pi\in\Pi$ there exists a policy $\pi'\in \Pi_A$ such that for any $x\notin A$ and $t\geq0$, it holds that $P_x^{\pi}(\tau_A = t) =  P_x^{\pi'}(\tau_A = t)$.
\end{lemma}
\begin{pf}
Let $\pi = \{\mu_t\}_{t\geq 0}\in\Pi$ and $\pi_A= \{\mu^A_{t}\}_{t\geq 0}\in \Pi_A$, which exists by Assumption \ref{assumA}. Define a policy $\pi'= \{\mu'_t\}_{t\geq 0}\in\Pi_A$ as follows:
\[\mu'_t(\cdot | h_t )=
\begin{cases}
 \mu_t (\cdot | h_t ), \mbox{ if } \overline{x}^t \notin A \\
 \mu^A_{t}(\cdot | h_t ), \mbox{ otherwise.}
\end{cases}
\]
Clearly $\pi'\in\Pi_A$. Let $x\notin A$, note that
\begin{align*}
 P_x^{\pi}&(\tau_A = t)  = P_x^{\pi}(X_1\notin A,\ldots ,X_{t-1}\notin A, X_{t}\in A) \\
   &= \sum_{\overline{x}^{t-1} \notin A,x_{t}\in A} P_{x}^{\pi}(X_1=x_1,\ldots , X_{t}=x_{t}) \\
    &= \sum_{\overline{x}^{t-1} \notin A,x_{t}\in A, u_i\in\mathbb{U}(x_i)}  P_{x}^{\pi}(u_0, x_1,u_1 \ldots x_{t},u_t) \\
   &=  \sum_{\overline{x}^{t-1} \notin A,x_{t}\in A, u_i\in\mathbb{U}(x_i)} P^{\pi}(x_t,u_t | h_{t-1}, u_{t-1}) \cdots \\
&\quad\quad\quad P^{\pi}(x_1,u_1 | h_{0}, u_{0} ) P^{\pi}(u_0 | x ).
\end{align*}
Since $x\notin A$ we have that $P^{\pi}(u_0 | x ) = \mu_{0}(u_0|x) = \mu'_{0}(u_0|x) = P^{\pi'}(u_0 | x )$. Also, for all $i\geq 1$,
\begin{align*}
  P^{\pi}(x_i,u_i | h_{i-1}, u_{i-1}) & = P^{\pi}(u_i | h_{i})P^{\pi}(x_i| h_{i-1}, u_{i-1})\\
  & = \mu_{i}(u_i|h_i)Q(x_i|x_{i-1},u_{i-1}).
\end{align*}
Because of the definition of $\pi'$ we know that $\mu_{i}(u_i|h_i) = \mu'_{i}(u_i|h_i)$ whenever $\overline{x}^t \notin A $. Thus, the fact that $x,x_1\ldots x_{t-1}\notin A$ and previous equality imply that for $i\leq t$ we have that $P^{\pi}(x_i,u_i | h_{i-1}, u_{i-1})  = P^{\pi'}(x_i,u_i | h_{i-1}, u_{i-1}).$ As a consequence $P_x^{\pi}(\tau_A = t) =  P_x^{\pi'}(\tau_A = t).$
\end{pf}

To proof the proposition let $\pi'$ be defined as in Lemma~\ref{NewPolicy}. If $x \in A$, the result follows trivially since $P_x^{\pi'}(\tau_A \leq t)=1$ for all $t\geq 0$. Let $x\notin A$. By of Lemma~\ref{NewPolicy} we know that $P_x^{\pi}(\tau_A \leq t)=P_x^{\hat{\pi}}(\tau_A \leq t),$ and since $\{X_t\in A \} \subset \{ \tau_A\leq t\}$, we obtain that $P_x^{\pi}(X_t\in A) \leq P_x^{\pi}(\tau_A \leq t) =  P_x^{\pi'}(\tau_A \leq t).$ Therefore, taking $\liminf$ and using Proposition \ref{HtA}, we get
\begin{align*}
\liminf_{t\rightarrow \infty}  P_x^{\pi}(X_t\in A) &\leq  \lim_{t\rightarrow \infty}  P_x^{\pi'}(X_t\in A) \\
&=P_x^{\pi'}(\tau_A<\infty ).
\end{align*}

\section{Proof of Proposition \ref{propequi}}\label{appB}

Equation~\eqref{equivalent} follows from combining the next two lemmas.
\begin{lemma}\label{1}
Given a policy $\pi\in \Pi$, $P^{\pi}_{\nu}(\tau_A<\tau_B,\tau_A<\infty)=\widetilde{P}^{\pi}_{\nu}(\tau_A<\tau_B,\tau_A<\infty).$
\end{lemma}
\begin{pf} If $x\in A$ both probabilities are 1 and if $x\in B$ both are 0. Let $x \in \mathbb{X}\setminus\{A \cup B\}$. It's clear that, $\{\tau_A<\tau_B ,\tau_A<\infty\}=\bigcup_{t\in\mathbb{N}} \{\tau_A=t, t<\tau_B\}$. Moreover, the events $\{\tau_A=t, t<\tau_B\}$ are disjoint, therefore to obtain the desired equality it is enough to show that $P_x^{\pi}(\tau_A = t, t < \tau_B) =  \widetilde{P}^{\pi}_{x}(\tau_A = t, t < \tau_B).$ Note that, 
\begin{align*}
 P_x^{\pi}&(\tau_A = t, t < \tau_B) \\
& = P_x^{\pi}(X_1\notin A\cup B,\ldots ,X_{t-1}\notin A\cup B, X_{t}\in A) \\
   &= \sum_{\overline{x}^{t-1} \notin A\cup B,x_{t}\in A} P_{x}^{\pi}(X_1=x_1, \ldots,  X_{t}=x_{t}) \\
       &= \sum_{\overline{x}^{t-1} \notin A\cup B,x_{t}\in A, u_i\in\mathbb{U}(x_i)}  P_{x}^{\pi}(u_0, x_1,u_1 \ldots x_{t},u_t) \\
   &=  \sum_{\overline{x}^{t-1} \notin A\cup B,x_{t}\in A, u_i\in\mathbb{U}(x_i)} P^{\pi}(x_t,u_t | h_{t-1}, u_{t-1}) \cdots\\ &\quad\quad\quad\quad P^{\pi}(x_1,u_1 | h_{0}, u_{0} ) P^{\pi}(u_0 | x ).
\end{align*}
We have $P^{\pi}(u_0 | x ) = \mu_{0}(u_0|x) = \widetilde{P}^{\pi}(u_0 | x )$. Furthermore, for all $i\geq 1$
\begin{align*}
  P^{\pi}(x_i,u_i | h_{i-1}, u_{i-1}) & = P^{\pi}(u_i | h_{i}) P^{\pi}(x_i| h_{i-1}, u_{i-1})\\
  & = \mu_{i}(u_i|h_i)Q(x_i|x_{i-1},u_{i-1}).
\end{align*}
Because of the definition of $\widetilde{Q}$ we know that $Q(\cdot|s,u) =\widetilde{Q}(\cdot|s,u)$ whenever $s\notin A\cup B$. Thus, using the previous equality and the fact that $x,x_1\ldots x_{t-1}\notin A\cup B$ we obtain that for $i\leq t$, $P^{\pi}(x_i,u_i | h_{i-1}, u_{i-1})  = \widetilde{P}^{\pi}(x_i,u_i | h_{i-1}, u_{i-1}),$ and therefore $P_x^{\pi}(\tau_A = t, t < \tau_B) =  \widetilde{P}^{\pi}_{x}(\tau_A = t, t < \tau_B).$
\end{pf}

\begin{lemma}\label{2}
Given a policy $\pi\in\Pi$, it holds that $\widetilde{P}^{\pi}_{\nu}(\tau_A<\tau_B , \tau_A<\infty)=\widetilde{P}^{\pi}_{\nu}(\tau_A<\infty).$
\end{lemma} 
\begin{pf}
Consider the event $\{\tau_A<\infty,\tau_B<\infty\}$. Let us prove that $\widetilde{P}^{\pi}_{x}(\tau_A<\infty ,\tau_B<\infty)=0$ for any $x\in\SX$. First of all, it is clear that $\{\tau_A<\infty ,\tau_B<\infty\}=\bigcup_{t,s\in\mathbb{N}} \{\tau_A=t, \tau_B=s\}.$ Note that for all $s,t$ the events $\{\tau_A=t, \tau_B=s\}$ are disjoint and are contained in either one of the following sets:
\begin{align*}
\{\exists n,m \in \mathbb{N} \mbox{ such that } n\leq m, X_n\in A, X_m\notin A\},\\
\{\exists n,m \in \mathbb{N} \mbox{ such that } n\leq m, X_n\in B, X_m\notin B\}.
\end{align*}
By the definition of $\widetilde{Q}$, the sets $A,B$ are closed under $\pi$. Therefore, Lemma~\ref{nullevent} implies that $\widetilde{P}^{\pi}_{x}(\tau_A=t, \tau_B=s)=0$, hence $\widetilde{P}^{\pi}_{x}(\tau_A<\infty ,\tau_B<\infty)=0$. Now, note that
\begin{align*}
\{\tau_A<\tau_B , \tau_A<\infty\}=&\{\tau_A<\tau_B , \tau_A<\infty,\tau_B<\infty\}\\
&\cup \{\tau_A<\tau_B , \tau_A<\infty\ , \tau_B=\infty\}\\
= &\{\tau_A<\tau_B , \tau_A<\infty,\tau_B<\infty\}\\
&\cup \{\tau_A<\infty\ , \tau_B=\infty\}.
\end{align*}
Similarly, $\{\tau_A<\infty\}=\{\tau_A<\infty,\tau_B<\infty\}\cup \{\tau_A<\infty\ , \tau_B=\infty\}.$ Therefore,
\begin{align*}
\widetilde{P}^{\pi}_{x}(\tau_A<\tau_B , \tau_A<\infty)&=\widetilde{P}^{\pi}_{x}(\tau_A<\infty , \tau_B=\infty)\\
&=\widetilde{P}^{\pi}_{x}(\tau_A<\infty).
\end{align*}
\end{pf}

\section{Proof of Lemma \ref{LemmaRestr}}\label{appLemma}

First, note that
\begin{align*}
P^{\pi}_x(h_t, u_t) & = P^{\pi}(x_t, u_t |h_{t-1}, u_{t-1} )\cdots \\
&\quad\quad\quad P^{\pi}(x_{1}, u_{1} |x, u_{0} ) P^{\pi}(u_0|x),
\end{align*}
and for $0<r\leq t$
\begin{align}\notag
 P^{\pi}(x_r,u_r | h_{r-1}, u_{r-1}) & = P^{\pi}(u_r | h_{r})P^{\pi}(x_r| h_{r-1}, u_{r-1})\\
  & = \pi(u_r|h_r)Q(x_r|x_{r-1},u_{r-1}). \label{E1}
\end{align}
\begin{enumerate}
\item[\underline{Case 1}:] Suppose $\overline{x}^t\notin B$. The fact that $x\notin B$ implies that $\pi(u_0|x) = \widehat{\pi}(u_0|(x,0))$, so that $P^{\pi}(u_0 | x) = \widehat{P}^{\widehat{\pi}}(u_0 | (x,0) )$.
Using equation \eqref{E1} and the definitions of $\widehat{Q}$ and $\widehat{\pi}$, we obtain that
\begin{align*}
 P^{\pi}&(x_r,u_r | h_{r-1}, u_{r-1}) = \pi(u_r|h_r)Q(x_r|x_{r-1},u_{r-1})\\
 & = \widehat{\pi}(u_r|\widehat{h}_r)\widehat{Q}((x_r,0)|(x_{r-1},0),u_{r-1})\\
 & = \widehat{P}^{\widehat{\pi}}(u_r | \widehat{h}_r)\widehat{P}^{\widehat{\pi}}((x_r,0)| \widehat{h}_{r-1}, u_{r-1})\\
 & = \widehat{P}^{\widehat{\pi}}((x_r,0),u_r| \widehat{h}_{r-1}, u_{r-1}).
\end{align*}
Hence $P^{\pi}_x(h_t, u_t)=\widehat{P}^{\widehat{\pi}}_{(x,0)}(\widehat{h_t}, u_t)$. 
\item[\underline{Case 2}:] Suppose $x_s\in B$ for some $0\leq s\leq t$, and $\overline{x}^{s-1}\notin B$. If $s=0$ then $x\in B$  and so $\pi(u_0|x) = \widehat{\pi}(u_0|(x,1))$ and $P^{\pi}(u_0 | x ) = \widehat{P}^{\widehat{\pi}}(u_0 | (x,1) )$. Otherwise, if $s>0$ we have that $P^{\pi}(u_0 | x ) = \widehat{P}^{\widehat{\pi}}(u_0 | (x,0) )$,  and as in Case 1, for $0<r<s$, we obtain that $P^{\pi}(x_r,u_r| h_{r-1}, u_{r-1}) = \widehat{P}^{\widehat{\pi}}((x_r,0),u_r | \widehat{h}_{r-1}, u_{r-1}).$ Now, for $r\geq s$, similarly, we obtain that 
$$P^{\pi}(x_r,u_r | h_{r-1}, u_{r-1})=\widehat{P}^{\widehat{\pi}}((x_r,1),u_r| \widehat{h}_{r-1}, u_{r-1}).$$ 
Therefore, $P^{\pi}_x(h_t, u_t)$ equals $\widehat{P}^{\widehat{\pi}}_{(x,0)}(\widehat{h_t}, u_t)$ if $x\notin B$, and  it equals $\widehat{P}^{\widehat{\pi}}_{(x,1)}(\widehat{h_t}, u_t)$ if $x\in B$. 
\end{enumerate}

\section{Proof of Theorem \ref{Linaverage}}\label{appC}

The proof of this theorem is divided into several lemmas. The first one allows us to rewrite the average reward function \eqref{avcost} in terms of the probabilities of hitting times being finite. Its proof is analogous to the proof of Theorem \ref{eqAv-Ht}.
\begin{lemma}
Let $\widehat{\pi}$ be a policy over $\widehat{\mathbb{X}}$. Then
$$
  v^{\widehat{\pi}}(x,i) = \widehat{P}_{(x,i)}^{\widehat{\pi}}\left(\tau_{ A\times\{ 0,1 \} }<\infty\right)
  - \lambda \widehat{P}_{(x,i)}^{\widehat{\pi}}\left( \tau_{ \SX\times\{1\} }<\infty\right). $$
\end{lemma}

In the following lemmas we will always assume that $\widehat{\pi}\in\widehat{\Pi}$ and that $\pi\in\Pi$ are given by Definition \ref{indpi}.

\begin{lemma}
If $x\notin B$, then $\widehat{P}_{(x,0)}^{\widehat{\pi}}( \tau_{ \SX\times\{ 1 \} }<\infty) = P_{x}^{\pi}( \tau_{ B } < \infty ).$
\end{lemma}
\begin{pf}
Note that the event $\{  \tau_{ B }<\infty  \}$ can be partitioned in disjoint events $\{  \tau_{ B}<\infty  \} = \bigcup_{s\in\mathbb{N}} \{  \tau_{ B} = s  \}.$ Therefore, it is enough to restrict our attention to such events. Using Lemma \ref{LemmaRestr} we obtain that
\begin{align}\notag
 &P_{x}^{\pi}( \tau_{B} = s  ) = \sum_{\overline{x}^{s-1}\notin B\atop x_{s}\in B, u_i\in\mathbb{U}(x_i)} P_{x}^{\pi}(u_0, x_1, u_1, \ldots, x_{s}, u_s)\\
 & = \sum_{{\overline{x}^{s-1}\notin B\atop x_{s}\in B}\atop u_i\in\mathbb{U}(x_i)}  \widehat{P}^{\widehat{\pi}}_{(x,0)}(u_0, (x_{1},0), u_1, \ldots,  (x_{s-1},0), u_{s-1}, (x_s,1), u_s). \label{E2}
\end{align}
Now,
\begin{align*}
\widehat{P}^{\widehat{\pi}}_{(x,0)} & (u_0, (x_{1},0), u_1, \ldots,  (x_{s-1},0), u_{s-1}, (x_s,1), u_s) \\ 
& = \widehat{P}^{\widehat{\pi}}((x_s,1),u_s|\widehat{h}_{s-1},u_{s-1}) \cdots \\
&\quad\quad\widehat{P}^{ \widehat{\pi}}((x_1,0),u_1|\widehat{h}_0,u_0) \widehat{P}^{ \widehat{\pi}}(u_0|(x,0)).
\end{align*}
Because of \eqref{E1} (which is valid for any Markov model) we have that, 
\begin{align*}
 \widehat{P}^{\widehat{\pi}}((x_t,i_t),u_t &| \widehat{h}_{t-1}, u_{t-1}) =\\
& \widehat{\pi}(u_t|\widehat{h_t})\widehat{Q}((x_t,i_t)|(x_{t-1},0),u_{t-1}).
 \end{align*}
If $x_t\in B$, for $t<s$ or if $x_s\notin B$, then the definition of $\widehat{Q}$ implies that $\widehat{P}^{\widehat{\pi}}((x_t,0),u_t | \widehat{h}_{t-1}, u_{t-1})=\widehat{P}^{\widehat{\pi}}((x_s,1),u_s | \widehat{h}_{s-1}, u_{s-1})=0$. Therefore, \eqref{E2} can be written as
\begin{align*}
 \sum_{ x_i\in \mathbb{X}\atop u_i\in\mathbb{U}(x_i)}  \widehat{P}^{\widehat{\pi}}_{(x,0)}(u_0, (x_{1},0), u_1, \ldots,&  (x_{s-1},0), u_{s-1}, (x_s,1), u_s)\\
& =  \widehat{P}^{\widehat{\pi}}_{(x,0)} ( \tau_{ \mathbb{X}\times\{1\} } = s).
\end{align*}
\end{pf}

Next two lemmas prove that whenever $x\notin B$, we have that $P_x^{\pi}(\tau_A<\infty) = \widehat{P}_{(x,0)}^{\widehat{\pi}}(\tau_{ A\times\{ 0,1 \} }<\infty).$

\begin{lemma}
If $x\notin B$, then  $\widehat{P}_{(x,0)}^{\widehat{\pi}}( \tau_{ A\times\{ 1 \} }<\infty) = P_{x}^{\pi}( \tau_{ A } < \infty , \tau_B<\tau_A).$
\end{lemma}
\begin{pf}
It is enough to prove that $P_{x}^{\pi}(\tau_{A} = t , \tau_B < \tau_A ) = \widehat{P}_{(x,0)}^{\widehat{\pi}}( \tau_{ A\times\{ 1 \} }=t)$. First, note that
 \begin{align*}
&  \widehat{P}_{(x,0)}^{ \widehat{\pi} } ( \tau_{ A\times \{1 \}}  = t ) \\
   &= \sum_{ u_i\in\mathbb{U}(j_i),k<t\atop j_k\notin A\times \{1 \}, j_t \in A \times \{1 \} } \widehat{P}_{(x,0),}^{ \widehat{\pi} } (u_0, j_1, u_1 \ldots,j_t, u_t)\\ 
&=\sum_{ u_i\in\mathbb{U}(j_i),k<t\atop j_k\notin A\times \{1 \}, j_t \in A \times \{1 \} } \widehat{P}^{\widehat{\pi}}(j_t,u_t|\widehat{h}_{t-1},u_{t-1}) \cdots \\
&\quad\quad\quad\quad \widehat{P}^{ \widehat{\pi}}(j_{1},u_1|\widehat{h}_0,u_0) \widehat{P}^{ \widehat{\pi}}(u_0|(x,0)).
 \end{align*}
Since the sets $\mathbb{X}\times \{1\}$ and $A\times\{0,1\}$ are closed for any policy $\widehat{\pi}$ by definition of $\widehat{Q}$, it is enough to consider histories of the form,
\[
((x,0),(x_1,0),\ldots,(x_{s-1},0),(x_{s},1), \ldots, (x_{t},1) )
\]
where $x_1,\ldots, x_{s-1}\notin B\cup A,x_s\in B, x_{s+1},\ldots, x_{t-1}\notin A, x_t\in A$, so by Lemma \ref{LemmaRestr} the last summation above is equal to
\begin{align*}
  \sum_{ { 1\leq s<t\atop \overline{x}^{s-1}\notin B\cup A, u_i\in\mathbb{U}(x_i)}\atop x_{s}\in B,x_{s+1},\ldots , x_{t-1}\notin A , x_{t}\in A } & P_{x}^{\pi}(u_0,x_1,\ldots,x_{t}, u_t)\\
&= \sum_{ s<t} P_{x}^{\pi}( \tau_{B} = s ,\tau_A=t ) \\ 
&=P_{x}^{\pi}(\tau_{A} = t , \tau_B < \tau_A ).
 \end{align*}
\end{pf}
Similarly, we obtain the next result.
\begin{lemma}
If $x\notin B$, then  $\widehat{P}_{(x,0)}^{\widehat{\pi}}( \tau_{ A\times\{ 0 \} }<\infty) = P_{x}^{\pi}( \tau_{ A } < \infty , \tau_A<\tau_B).$
\end{lemma}
The final lemma completes the proof of the theorem by considering the case when the initial state belongs to $B$.
\begin{lemma}
If $x\in B$, then $v^{\widehat{\pi}}(x,1) = P_{x}^{\pi}(\tau_A<\infty ) - \lambda.$
\end{lemma}
\begin{pf}
Clearly, it is enough to show that $P_{x}^{\pi}(\tau_A<\infty ) = \widehat{P}_{(x,1)}^{\widehat{\pi}}(\tau_{A\times \{1\} } <\infty )$. As before, consider the events $\{ \tau_{A} = t \}$. Since $x\in B$, Lemma \ref{LemmaRestr} implies that
\begin{align*}
 P_{x}^{\pi}&(\tau_A =t  ) = \sum_{\overline{x}^{t-1}\notin A, x_{t}\in A, u_i\mathbb{U}(x_i) } P^{\pi}_x(u_0, x_1, \ldots, x_t, u_t)\\
 & = \sum_{\overline{x}^{t-1}\notin A, x_{t}\in A, u_i\mathbb{U}(x_i) } \widehat{P}^{\widehat{\pi}}_{(x,1)}(u_0, (x_1,1), \ldots, (x_t,1), u_t)\\
&=\sum_{u_i\in\mathbb{U}(j_i),k<t\atop j_k\notin A\times \{1 \}, j_t \in A \times \{1 \} }  \widehat{P}^{\widehat{\pi}}_{(x,1)}(u_0, j_1, \ldots, j_t, u_t)\\
 &=  \widehat{P}_{(x,1)}^{ \widehat{\pi} }  ( \tau_{ A\times \{1 \}}  = t ),
\end{align*}
where the second to the last equality follows from the fact that $\mathbb{X}\times \{ 1 \}$ is closed.
\end{pf}

\end{document}